\documentclass{cup-elements_dat}
\usepackage{blindtext}
  \usepackage{amsthm}
 \usepackage{url}
\usepackage{makeidx}
\makeindex

\usepackage{newtxtext}
 \usepackage[varbb]{newtxmath}
\usepackage[scaled=0.86]{helvet}
\usepackage{multirow}
\usepackage{graphicx}

\usepackage{natbib}
\usepackage[caption=false]{subfig}
\usepackage{amsthm}
\usepackage{amsmath}
\usepackage{amsmath,bm}
\usepackage{amssymb}
\usepackage{amsfonts}
\usepackage[export]{adjustbox}
\usepackage{braket,amsfonts}
\usepackage{amsmath}
\usepackage{enumitem}
\newcommand*\eval[3]{\left.#1\right\rvert_{#2}^{#3}}
\usepackage{mathtools}
\usepackage{mismath}
\DeclarePairedDelimiterX{\inp}[2]{\langle}{\rangle}{#1, #2} 
\usepackage{algorithm,algpseudocode}
\algnewcommand{\Inputs}[1]{%
  \State \textbf{Inputs:}
  \Statex \hspace*{\algorithmicindent}\parbox[t]{.8\linewidth}{\raggedright #1}
}
\algnewcommand{\Initialize}[1]{%
  \State \textbf{Initialize:}
  \Statex \hspace*{\algorithmicindent}\parbox[t]{.8\linewidth}{\raggedright #1}
} 
\algnewcommand{\Outputs}[1]{%
  \State \textbf{Outputs:}
  \Statex \hspace*{\algorithmicindent}\parbox[t]{.8\linewidth}{\raggedright #1}
}
\usepackage{xcolor}

\usepackage{acronym}
\acrodef{FTFC}[FTFC]{\emph{Fundamental Theorem of Fractional Calculus}}
\acrodef{KMD}[KMD]{\emph{Koopman Mode Decomposition}}
\acrodef{KEF}[KEF]{\emph{Koopman Eigenfunction}}
\acrodef{KEFal}[KEFal]{\emph{Koopman Eigenfunctional}}
\acrodef{DMD}[DMD]{\emph{Dynamic Mode Decomposition}}
\acrodef{SDMD}[S-DMD]{\emph{Symmetric DMD}}
\acrodef{EDMD}[EDMD]{\emph{Extended DMD}}
\acrodef{KDMD}[KDMD]{\emph{kernel DMD}}
\acrodef{SVD}[SVD]{\emph{Singular Value Decomposition}}
\acrodef{TV}[TV]{\emph{Total Variation}}
\acrodef{DFT}[DFT]{\emph{Discrete Fourier Transform}}
\acrodef{CFT}[CFT]{\emph{Continuous Fourier Transform}}
\acrodef{POD}[POD]{\emph{Proper Orthogonal Decomposition}}
\acrodef{ROA}[ROA]{\emph{Region of Attraction}}
\acrodef{PDE}[PDE]{\emph{Partial Differential Equation}}
\acrodef{PDEs}[PDEs]{\emph{Partial Differential Equations}}
\acrodef{TV}[TV]{\emph{Total Variation}}
\acrodef{OrthoNS}[OrthoNS]{\emph{Orthogonal Nonlinear Spectral decomposition}}

\usepackage{longtable}

\CUPseries{Ido Cohen \& Guy Gilboa}
\CUPelements{Latent Modes of Nonlinear Flows - a Koopman Theory Analysis}

\title{Latent Modes of Nonlinear Flows - a Koopman Theory Analysis}

\author{Ido Cohen}
\affil{Faculty of Electrical And Computer Engineering - Technion}

\author{Guy Gilboa}
\affil{Faculty of Electrical And Computer Engineering - Technion}

  \theoremstyle{plain}
  \newtheorem{theorem}{Theorem}[section]
\newtheorem{lemma}[theorem]{Lemma}
  \newtheorem{proposition}[theorem]{Proposition}
  
  \newtheorem{corollary}[theorem]{Corollary}
  \newtheorem*{theorem*}{Theorem}
  \newtheorem*{lemma*}{Lemma}
  \newtheorem*{proposition*}{Proposition}
  \newtheorem*{corollary*}{Corollary}
  \newtheorem*{conjecture*}{Conjecture}

  \theoremstyle{definition}
  \newtheorem{definition}[theorem]{Definition}
  
  \newtheorem{remark}[theorem]{Remark}

\newtheorem{idoexample}[theorem]{Example}
\newtheorem{assumption}[theorem]{Assumption}
\begin{document}

\frontmatter  
\maketitle

\begin{abstract}
Extracting the latent underlying structures of complex nonlinear local and nonlocal flows is essential for their analysis and modeling. 
In this work we attempt to provide a consistent framework through Koopman theory and its related popular discrete approximation -- dynamic mode decomposition (DMD). We investigate the conditions to perform appropriate linearization, dimensionality reduction and representation of flows in a highly general setting.

The essential elements of this framework are \ac{KEF}, for which existence conditions are formulated. This is done by viewing the dynamic as a curve in state-space. These conditions lay the foundations for system reconstruction, global controllability, and observability for nonlinear dynamics.

We examine the limitations of DMD through the analysis of Koopman theory and propose a new mode decomposition technique based on the typical time profile of the dynamics. An overcomplete dictionary of decay profiles is used to sparsely approximate the flow. This analysis is also valid in the full continuous setting of Koopman theory, which is based on variational calculus.
We demonstrate applications of this analysis, such as finding \acp{KEF} and their multiplicities, dynamics reconstruction and global linearization.

\end{abstract}



\copyrightauthor{Ido Cohen, Guy Gilboa, 2021}

\mainmatter  

{\bf{List of abbreviations}}
\addcontentsline{toc}{section}{List of abbreviations}
\begin{longtable}{lp{0.6\textwidth}}
    \acs{DMD}&\acl{DMD}\\
    \acs{EDMD}&\acl{EDMD}\\
    \acs{KDMD}&\acl{KDMD}\\
    \acs{KEF}&\acl{KEF}\\
    \acs{KEFal}&\acl{KEFal}\\
    \acs{KMD}&\acl{KMD}\\
    \acs{PDE}&\acl{PDE}\\
    \acs{ROA}&\acl{ROA}\\
    \acs{SDMD}&\acl{SDMD}\\
    \acs{SVD}&\acl{SVD}\\
    \acs{TV}&\acl{TV}\\

  \end{longtable}
\section{Introduction}
Knowing the latent space of certain data allows one to represent it concisely and to differentiate between signal and clutter parts. Recovering this space in a data-driven manner is a long-standing research problem. Data resulting from dynamical systems is represented commonly as spatial structures (modes) that are attenuated or enhanced with time. 
A common technique in linear flows is \emph{separation of variables}. It is assumed that a solution $u(x,t)$ of a linear flow can be expressed as,  
\begin{equation}
    u(x,t) = X(x)T(t).
\end{equation}
That is, the solution is a multiplication of a function of the spatial variable $x$ and a function of the temporal variable $t$. 
In this study we examine, from various angles, the following paradigm: a \emph{nonlinear flow} can be well approximated (or even exactly expressed) by a linear combination of variable separated functions, 
\begin{equation}
    u(x,t) \approx \sum_{i=1}^m X_i(x)T_i(t).
\end{equation}
In this context, the spatial structures $X_i$ are referred to as \emph{modes} and  $T_i$ are \emph{time-profiles}. For such an approximation, if the error is negligible and $m$ is small, we obtain a significant simplification of the system. This enables better understanding and modeling, allowing accurate interpolation and prediction of the dynamics.  

The theory of Koopman argues that for many nonlinear systems  data measurements evolve as if the dynamical system is linear (in some  infinite dimensional space). A well-known algorithm to approximate these measurements is \acf{DMD} of \cite{schmid2010dynamic}. In this work, we formulate sufficient and necessary conditions for the existence of these measurements. These findings highlight certain flaws of \ac{DMD}. Finally, we suggest a new mode decomposition to overcome some of these problems, originated in an algorithm for general spectral decomposition of \cite{gilboa2018nonlinear}.

In many dynamical processes, there are measurements of the observations that evolve linearly, or approximately so, see \cite{otto2021koopman}. A theoretical justification for that can be traced back to the seminal work of \cite{koopman1931hamiltonian}. These measurements are referred to as \acfp{KEF}. An algorithm was proposed by \cite{mezic2005spectral}, \acf{KMD}, to reconstruct the dynamics using spatial structures, termed as modes, which are the coefficients of Koopman eigenfunctions. Since \acp{KEF} evolve as if they were observations in a linear dynamical system, \ac{KMD} can interpret the original dynamics as a linear one.

This decomposition might be infinite-dimensional. In \cite{schmid2010dynamic} \ac{DMD} it was suggested to approximate \ac{KMD} in a finite domain. If the \acp{KEF} measurements are linear combinations of the observations then \ac{DMD} yields the Koopman mode accurately. As noted in \cite{kutz2016dynamic}, \ac{DMD} can be interpreted as an exponential data fitting algorithm.  In the more general nonlinear case, \ac{DMD} may not reveal well the underlying modes and the dynamics.

Recently the authors and colleagues have formalized this insight in \cite{cohen2021modes}, in the context of homogeneous flows, referring to it as the \ac{DMD} paradox. As the step-size approaches zero, dynamic reconstruction with \ac{DMD} results in positive mean squared error, but, paradoxically, with zero \ac{DMD} error. In general, this paradox exists in any dynamical system whose \acp{KEF} are not linear combinations of the observations. This phenomenon becomes extreme when the system is zero homogeneous, as shown in \cite{cohen2021Total}. Such cases are common in gradient flows of one-homogeneous functionals, such as local or nonlocal TV-flows, \cite{andreu2001minimizing}, \cite{gilboa2009nonlocal}. In that case, the dynamics is only in $C^0$ almost everywhere and exponential decay is a very crude and inaccurate approximation. For such flows, lifting the observations to a finite higher dimensional space does not solve the problem (see for example \cite{korda2018linear,williams2015data2}).

This alleged contradiction between \ac{KMD} and \ac{DMD} leads us to examine the fundamentals of Koopman theory. We follow the general solution of a \ac{KEF} with respect to time and analyze the mapping between the state-space and the time variable. The existence of this mapping depends on the smoothness properties of the dynamics.

As a direct result, we introduce a new method that overcomes the \ac{DMD} limitations for smoothing-type processes. These findings, with some adaptations, are valid in the full continuous settings, as discussed by \cite{kutz2016koopman,mauroy2021koopman}.

\paragraph{Main Contributions}
We formulate the conditions for the existence of a \ac{KEF}. If it exists, there is an infinite set of \acp{KEF}. We distinguish between different types of eigenfunction groups and analyze their multiplicity. We show that certain multiplicities are crucial to obtain dynamics reconstruction, controllability, and observability (Section \ref{sec:disKoopman}). These conclusions are extended to the full continuous setting. Conditions for the existence of \acp{KEFal} are presented (Section \ref{sec:conKoopman}). Following these insights, we suggest an alternative algorithm for finding Koopman modes induced by fitting time profiles that best characterize the dynamics. This algorithm overcomes some inherent limitations of  \ac{DMD} (Section \ref{sec:ourDMD}). We attempt to bridge between nonlinear spectral decomposition and \ac{KMD}. Specifically, we show that spectral \ac{TV} of \cite{gilboa2014total} and its generalizations yield Koopman modes.
Throughout this work, we illustrate the theory with simple toy examples. Additional examples and experiments are given in Section \ref{sec:results} . In the following section, we provide the essential definitions and notations.

\section{Preliminary}
In this section, we present some background on Koopman operators, its eigenfunctions and eigenfunctionals and the related DMD framework. We note certain properties of variational calculus which are relevant to  Section \ref{sec:conKoopman}. In addition, we outline the work of \cite{gilboa2018nonlinear} and \cite{Katzir2017Thesis}, where nonlinear flows  are decomposed through a dictionary of decay profiles. We adapt this method for the extraction of Koompan modes in Section \ref{sec:ourDMD}.

\subsection{Koopman theory}
\subsubsection{Discrete spatial setting}
We consider a dynamical system in a semi-discrete setting, expressed as,
\begin{equation}\label{eq:disDynamicalSystem}
    \frac{d}{dt}{\bm{x}}(t)=P(\bm{x}(t)),\quad \bm{x}(0)=\bm{x}_0,\quad t\in I,
\end{equation}
where $\bm{x}\in \mathbb{R}^N$ is a state vector, $P:\mathbb{R}^N\to \mathbb{R}^N$ is a (nonlinear) operator, and $I=[a,b]\subseteq \mathbb{R}^+$ is a time interval.
Let $g:\mathbb{R}^N \to \mathbb{R}$ be a measurement of $\bm{x}$. The Koopman operator $K_P^\tau$ is a linear operator that acts on the infinite-dimensional space
of measurements $g(\bm{x})$ of the state, defined by \cite{koopman1931hamiltonian,mezic2005spectral},
\begin{equation}\label{eq:koopDef}\tag{\bf{KO}}
    K_P^\tau(g(\bm{x}(s)))=g(\bm{x}(s+\tau)),\quad s,s+\tau\in I,
\end{equation}
where $\tau>0$. The Koopman operator is linear, namely it admits,
\begin{equation*}
    K_P^\tau(\alpha g(\bm{x}(s))+\beta f(\bm{x}(s)))=\alpha K_P^\tau( g(\bm{x}(s)))+\beta K_P^\tau(f(\bm{x}(s))),
\end{equation*}
for all measurements $g$ and $f$ and for all constants $\alpha$ and $\beta$. In addition, the Koopman operators $\{K_P^\tau\}_{\tau\ge 0}$ admits a semigroup property, more formally,
\begin{equation*}
    K_P^{\tau_2}\circ K_P^{\tau_1}=K_P^{\tau_1+\tau_2},
\end{equation*}
where $\circ$ denotes the composition operator. 
An eigenfunction of the Koopman operator, $\varphi(\bm{x})$, admits,
\begin{equation}\label{eq:KEFdef}
    K_P^\tau(\varphi(\bm{x}(s)))=\varphi(\bm{x}(s+\tau))=\eta^\tau \varphi(\bm{x}(s)),\quad s,s+\tau\in I,
\end{equation}
for some $\eta\in \mathbb{C}$.
Due to the semigroup attribute of the Koopman operator, the following limit exists,
\begin{equation}
    \lim_{\tau\to 0}\frac{K_P^\tau(\varphi(\bm{x}(s)))-\varphi(\bm{x}(s))}{\tau}=\lim_{\tau\to 0}\frac{\varphi(\bm{x}(s+\tau))-\varphi(\bm{x}(s))}{\tau}=\eval{\frac{d}{dt}\varphi(\bm{x}(t))}{t=s}{}.
\end{equation}
This limit can be explained by the relations of the Koopman operator and Lie derivatives, see \cite{brunton2021modern}. It can be shown (see for instance  \cite{mauroy2020koopman}, p. 10) that a \ac{KEF} admits,
\begin{equation}\label{eq:KEFder}
    \frac{d}{dt}\varphi(\bm{x}(t)) = \lambda\cdot \varphi(\bm{x}(t)), \quad \forall t\in I,
\end{equation}
for some $\lambda\in \mathbb{C}$. The relation between $\eta$ in Eq. \eqref{eq:KEFdef} and $\lambda$ in Eq. \eqref{eq:KEFder} is $\eta = e^{\lambda}$.  The solution of this linear ODE is given by,
\begin{equation}\label{eq:KEF}\tag{KEF}
    \varphi(\bm{x}(t)) = \varphi(\bm{x}(a))e^{\lambda t} , \quad \forall t\in I.
\end{equation}

\paragraph{\acl{KMD}}
\acf{KMD} is a spatiotemporal mode decomposition of dynamical systems based on \acp{KEF}. Namely, the state space ${\bm{x}}$ can be expressed as (\cite{mezic2005spectral}),
\begin{equation}
	{\bm{x}}(t)=\sum_{i=1}^\infty\bm{v}_i\varphi_{i}(t),
\end{equation}
where $\varphi_i(t)$ is a \ac{KEF} and $\bm{v}_i$ is the corresponding vector, referred to as Koopman mode. When the dynamic is nonlinear the decomposition may have infinite elements. In practice, a finite approximation method is used.
The most common one is \ac{DMD}, as explained in Section \ref{sec:PreDMD}.

\subsubsection{Full continuous setting}
Let $u:L\subset\mathbb{R}\times I\subseteq\mathbb{R}^+$  
be the solution of the following PDE,
\begin{equation}\label{eq:PDE}
    u_t(x,t)=\mathcal{P}(u(x,t)), \quad u(x,0)=f(x).
\end{equation}
We assume that $u$ belongs to a Hilbert space $\mathcal{H}$ with an inner product, $\inp{v}{u}$ and its associated norm $\norm{\cdot}=\sqrt[]{\inp{\cdot}{\cdot}}$, $\mathcal{P}:\mathcal{H}\to \mathcal{H}$ is a (nonlinear) operator. Let $Q:\mathcal{H}\to \mathbb{R}$ be a proper, lower-semicontinuous functional.
The Koopman operator, $K_\mathcal{P}^t$, in the sense of PDE, is defined by \cite{nakao2020spectral},
\begin{equation}
    K_\mathcal{P}^\tau(Q(u(x,s)))=Q(u(x,s+\tau)), \quad s,s+\tau \in I.
\end{equation}
An eigenfunctional, $\phi$, of the Koopman operator is a functional admitting the following,
\begin{equation}\label{eq:KEFalDef}
    K_\mathcal{P}^\tau(\phi(u(x,s)))=\phi(u(x,s+\tau))=\eta^\tau \phi(u(x,s)), \quad s,s+\tau \in I.
\end{equation}
By letting $\tau \to 0$ an eigenfunctional of the Koopman operator admits the following ODE,
\begin{equation}\label{eq:eigenFunDer}
\frac{d}{dt}\phi(u(x,t))=\lambda \phi(u(x,t)),
\end{equation}
for some $\lambda\in \mathbb{C}$. The relation between $\eta$ in Eq. \eqref{eq:KEFalDef} and $\lambda$ in Eq. \eqref{eq:eigenFunDer} is $\eta = e^{\lambda}$. Thus, a \acf{KEFal} is of the form,
\begin{equation}\label{eq:KEFal}\tag{\bf{KEFal}}
    \phi(u(x,t)) = \phi(u(x,a))e^{\lambda t} , \quad \forall t\in I.
\end{equation}

\acl{KMD}. 
In the same manner as in the semi-discrete setting , we formulate the solution of the PDE, Eq. \eqref{eq:PDE}, with \acp{KEFal} (\cite{nakao2020spectral}). Namely, the solution $u(x,t)$ can be expressed as,
\begin{equation}
	u(x,t)=\sum_{i=1}^\infty d_i(x)\phi_{i}(u),
\end{equation}
where $\phi_i(u)$ is a \ac{KEFal} and $d_i(x)$ is the spatial mode. One way to approximate these spatial modes, is by the method introduced by \cite{nathan2018applied}.

\subsection{\acf{DMD}}\label{sec:PreDMD}
\ac{DMD} extracts the main spatial structures in the dynamics, \cite{schmid2010dynamic}. Backed by Koopman theory,
\ac{DMD} is a principal method to approximate the Koopman modes. It is a data driven method, based on snapshots (mostly, uniformly in time) of the dynamics, $\bm{x}_k=\bm{x}(t_k)$. The main steps in \ac{DMD} and its extensions (e.g. Exact \ac{DMD} \cite{tu2013dynamic}, tls\ac{DMD} \cite{hemati2017biasing}, fb\ac{DMD} \cite{dawson2016characterizing}, \acs{SDMD} \cite{cohen2021modes}, and optimized \ac{DMD}  \cite{askham2018variable}) are:
\begin{enumerate}
    \item \emph{Coordinates representation} - finding the main structures in the dynamics.
    \item \emph{Dimensionality reduction} - choosing the dominant parts of the dynamics.
    \item \emph{Linear mapping} - finding a linear mapping in the reduced dimensional space.
\end{enumerate}
We describe these steps in detail in Appendix \ref{appsec:DMDsteps}. The result of \ac{DMD} and its variants is sets of \emph{modes, $\{\bm{\phi}_i\}$, eigenvalues $\{\mu_i\}$, and coefficients $\{\alpha_i\}$}, where $i=1,..\,,r$ and $r$ is the reduced dimension. In the \ac{DMD} framework, the dynamics is approximated by,
\begin{equation}\label{eq:dynamicApproDiscrete}
    \bm{\tilde{x}}_k\approx \sum_{i=1}^r\alpha_i\mu_i^k\bm{\phi}_i.
\end{equation}

\subsection{General Spectral Decomposition} One of the main goals of signal analysis is to represent a signal sparsely,  yet precisely. We focus here on approximating a solution to a PDE, \eqref{eq:PDE}, by a decomposition of the form, 
\begin{equation}\label{eq:PDESolution}
    u(x,t)\approx \sum_{i=1}^{L}h_i(x)a_i(t),
\end{equation}
where $\{h_i(\cdot)\}_{i=1}^{L}$ are spatial functions and $\{a_i(t)\}_{i=1}^L$ are their respective time profiles. The time profiles are typical to the operator $\mathcal{P}$ and for homogeneous operators can be expressed analytically, see \cite{cohen:hal-01870019}.
In the semi-discrete setting, the approximate solution of Eq. \eqref{eq:disDynamicalSystem} can be expressed as,
\begin{equation}\label{eq:disFormSolDMD}
    \bm{x}(t)\approx \sum_{i=1}^{L}\bm{v}_ia_i(t),
\end{equation}
where $\{\bm{v}_i\}_{i=1}^L$ are spatial structures and $\{a_i(t)\}_{i=1}^L$ are the corresponding time profiles. Note that in some cases (e.g. linear diffusion or \ac{TV} flow, as shown in \cite{burger2016spectral}) Eqs. \eqref{eq:PDESolution} and \eqref{eq:disFormSolDMD} reach equality for finite or infinite $L$.

This is the basis of the general spectral decomposition suggested in the thesis of  \cite{Katzir2017Thesis} and summarized in the book of \cite{gilboa2018nonlinear} (chapter 9). The initial condition of Eq. \eqref{eq:PDE} is reconstructed with spatial structures that decay according to a known time profile. More formally, given the solution, $u(x,t)$, the spatial structures are the vectors of the minimizer of the following optimization problem,
\begin{equation}\label{eq:sparseKatzirEtGilboa}
    \min_{\mathcal{H}}\norm{\mathcal{U}-\mathcal{H}\mathcal{D}}_\mathcal{F}^2
\end{equation}
where $\mathcal{U}$ is a matrix of the sampled solution in  time and space, $\mathcal{H}$ is a matrix containing (in its columns) the main spatial structures, and $\mathcal{D}$ is a dictionary of decay profiles. One can formulate these matrices as
\begin{equation}\label{eq:Katzir}
\begin{split}
    \mathcal{U}&=\begin{bmatrix}u(x_1,t_0)&\cdots&u(x_1,t_M)\\
    \vdots&&\vdots\\
    u(x_N,t_0)&\cdots&u(x_N,t_M)\end{bmatrix},
    \mathcal{H}=\begin{bmatrix}h_1(x_1)&\cdots&h_r(x_1)\\
    \vdots&&\vdots\\
    h_1(x_N)&\cdots&h_r(x_N)\end{bmatrix},\\
    \mathcal{D}&=\begin{bmatrix}a_1(t_0)&\cdots&a_1(t_M)\\
    \vdots&&\vdots\\
    a_r(t_0)&\cdots&a_r(t_M)\end{bmatrix},
    \end{split}
\end{equation}
where $\mathcal{U}\in \mathbb{R}^{N\times (M+1)}$, $\mathcal{H}\in \mathbb{R}^{N\times r}$, and $\mathcal{D}\in \mathbb{R}^{r\times (M+1)}$.
The optimization problem, Eq. \eqref{eq:sparseKatzirEtGilboa}, fits also the form the of semi-discrete setting in the dynamics of Eq. \eqref{eq:disDynamicalSystem}, where it is sampled in the time axis. We can formulate the following optimization problem,
\begin{equation}
    \norm{X-\mathcal{V}\mathcal{D}}_\mathcal{F}^2,
\end{equation}
where the matrix $X$ contains the samples of the dynamics
\begin{equation}\label{eq:KatzirX}
    X=\begin{bmatrix}\bm{x}_0&\bm{x}_1&\cdots&\bm{x}_M\end{bmatrix}\in \mathbb{R}^{N\times (M+1)},
\end{equation}
the matrix $\mathcal{V}$ contains the main spatial structure of the dynamic (Eq. \eqref{eq:disFormSolDMD})
\begin{equation}\label{eq:KatzirV}
    \mathcal{V}=\begin{bmatrix}\bm{v}_1&\bm{v}_2&\cdots&\bm{v}_r\end{bmatrix}\in \mathbb{R}^{N\times r},
\end{equation}
and the dictionary, $\mathcal{D}$, remains unchanged.

\subsection{Variational Calculus}
\paragraph{Brezis chain rule}
Let $Q$ be a functional over some Banach space and $\partial Q$ be its variational derivative. Under the regime of the PDE, Eq. \eqref{eq:PDE}, we can formulate the time derivative of the functional, $Q(u(t))$, through the ``chain rule of Brezis'' \cite{brezis1973ope} as,
\begin{equation}
    \frac{d}{dt}Q(u(x,t))=\inp{\partial Q(u)}{\frac{d}{dt}u(x,t)}=\inp{\partial Q(u)}{\mathcal{P}(u(x,t))}.
\end{equation}

\paragraph{Fr\'echet Differentiability}
The operator $\mathcal{P}:\mathcal{H} \to \mathcal{H}$ is Fr\'echet differentiable at $u$ if there exists a bounded linear operator $\mathcal{L}$, such that,
\begin{equation}
    \lim_{\norm{h}\to 0}\frac{\norm{\mathcal{P}(u+h)-\mathcal{P}(u)-\mathcal{L}(h)}}{\norm{h}}=0
\end{equation}
holds from any $h\in \mathcal{H}$. In this case, $\mathcal{P}(u+h)$ can be expanded in the Landau notation as
\begin{equation}\label{eq:freDef}
    \mathcal{P}(u+h)=\mathcal{P}(u)+\mathcal{L}(h)+o(h),
\end{equation}
where $\lim_{\norm{h}\to 0}\norm{o(h)}/\norm{h}=0$. 

\paragraph{Proper Operator} The operator $\mathcal{P}(f(x))$ is proper if it gets a finite value for any $f(x)\in \mathcal{H}$ and for any $x\in[0,L]$.

\paragraph{\acf{ROA}} Let $\bm{x}^*$ be an equilibrium point of the dynamical system in Eq. \eqref{eq:disDynamicalSystem}. The region of attraction is the largest set in $\mathbb{R}^N$ that admits the following property: if the initial condition of the dynamics is from the set, then the system  converges to $\bm{x}^*$ (see e.g. \cite{valmorbida2017region}). More formally,
\begin{equation}\label{eq:ROA}
    \mathcal{RA}(\bm{x}^*)=\{x_{init}\in \mathbb{R}^N| \bm{x}(t=0)=x_{init}, \lim_{t\to \infty}\bm{x}(t)=\bm{x}^*\}.
\end{equation}

\section{Motivation for this work}

This monograph follows an earlier research, attempting to directly apply Koopman operator theory for homogeneous smoothing flows. 
In \cite{cohen2021modes} we investigated the use of \acf{DMD} for common nonlinear flows emerging in image processing, such as TV-flow and $p$-Laplacian flows. We found out that \ac{DMD} cannot be naively applied to decompose these flows and presented in detail certain flaws of this procedure.

\ac{DMD} has become a very common tool in dynamical system analysis. This decomposition provokes interest in many domains of research, such as  fluid dynamics, video processing, epidemiology, neuroscience, and finance, see \cite{kutz2016dynamic}. A main advantage is its simplicity and its ability to simplify complex processes by a few modes, in many cases. 

\ac{DMD} invokes well-established tools of dimensionality reduction, and can often reveal the main spatial components of the dynamic. However, the algorithm entails some fundamental problems in recovering nonlinear systems. Moreover, its drawbacks are emphasized when the dynamic is stable and we use the \ac{DMD} expansions such as, \cite{azencot2019consistent}, where the inverse dynamic is taken into account. Below we show some examples where \ac{DMD} is failing.

\subsection{\texorpdfstring{\ac{DMD}}{TEXT} paradox}
\begin{figure}[phtb!]
\centering
\includegraphics[width=0.95\textwidth,valign=c]{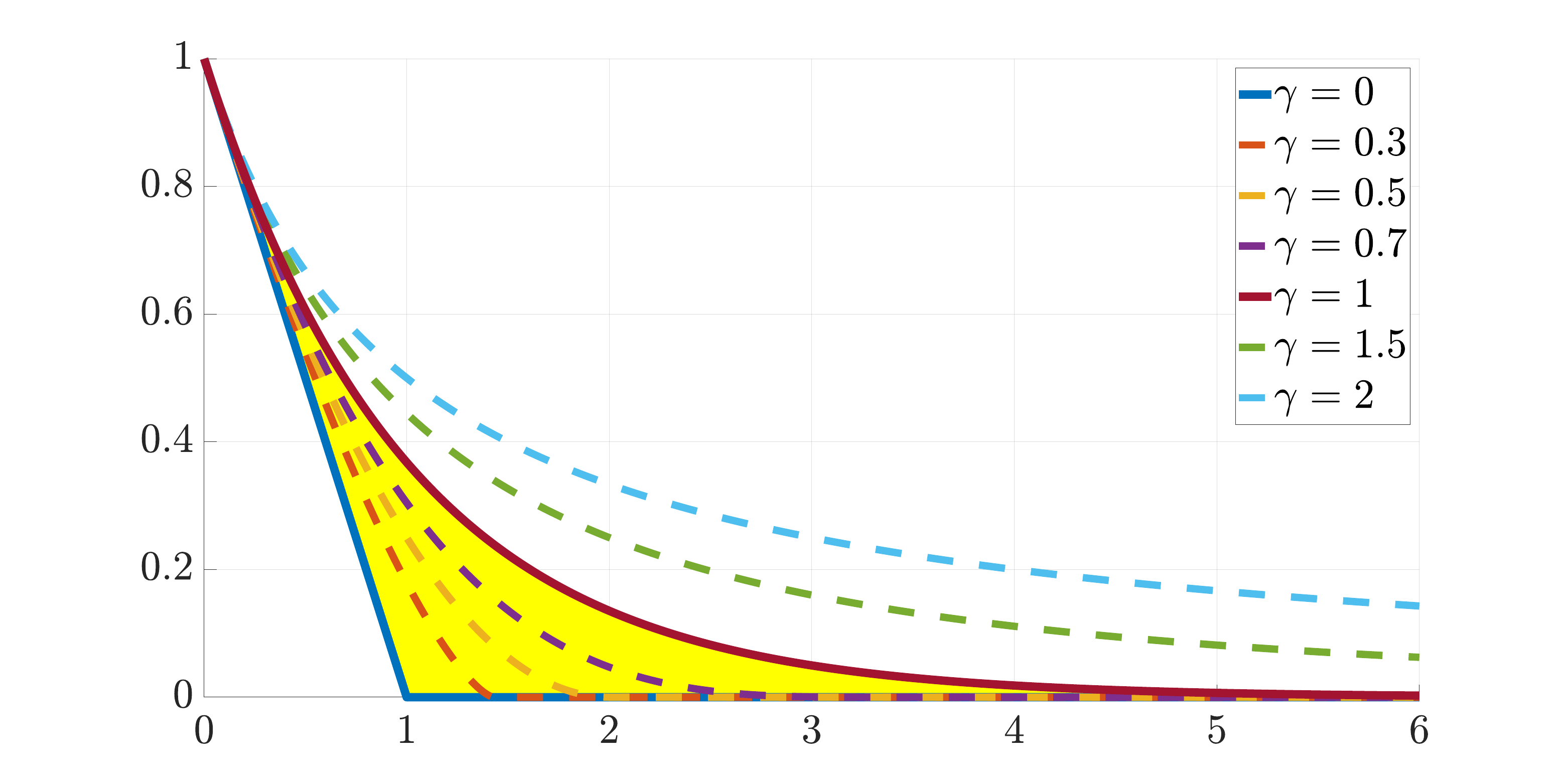}
\caption{{\bf{decay profile}}.   }
\label{Fig:decayProfile}
\end{figure}

The \ac{DMD} paradox was firstly introduced in \cite{cohen2021modes}. We recap here the findings about this \ac{DMD} flaw.
Let $P(\cdot)$ be a $\gamma$-homogeneous operator ($\gamma\in\mathbb{R}$) over some Banach space $\mathcal{B}$, i.e. $P(a v)=a\abs{a}^{\gamma -1}P(v)$ 
for any $a\in \mathbb{R}$ and $v\in\mathcal{B}$. Let $\phi\in \mathcal{B}$ be an eigenfunction of $P$, admitting $P(\phi)=\lambda \phi$ for a real valued $\lambda$. Then, the solution of the PDE (\cite{cohen:hal-01870019,cohen2020introducing})
\begin{equation}\label{motiEq:PDE}
    \frac{d}{dt}u=P(u), \quad u(t=0)=\phi,\, 
\end{equation}
is given by
\begin{equation}
u(t)=
\begin{cases}
    \phi\left[\left(1+\lambda (1-\gamma) t\right)^+\right]^{\frac{1}{1-\gamma}}& \gamma \ne 1\\
    \phi e^{\lambda t}& \gamma =1
\end{cases}
\end{equation}
where $(\cdot)^+=\max\{\cdot,0\}$. Under some conditions, the eigenvalue $\lambda$ is negative for any non-trivial eigenfunction $\phi$. Therefore, this solution gets the steady state in finite time when $\gamma \in [0,1)$. The time $T_{ext}$ for which the dynamic vanishes is 
\begin{equation}
    T_{ext}=\frac{1}{\lambda (\gamma -1)}.
\end{equation}
The decay profile is a fundamental characteristic of signal processing frameworks related to  eigenfunctions of $\gamma$-homogeneous operators, $\gamma\in[0,1)$. This decomposition generalizes the one based on gradient flows of one-homogeneous functionals, see \cite{bungert2019computing,burger2016spectral,cohen2020introducing,gilboa2013spectral,gilboa2014total}. The decay profile depends on the homogeneity order $\gamma$ (see Fig. \ref{Fig:decayProfile}). The decay varies from a truncated linear function for zero-homogeneous operators through truncated polynomial functions when $\gamma\in[0,1)$, to exponential function for one-homogeneous operators and finally to hyperbolic functions when $\gamma>2$.

The finite extinction time, inherent in flows where $\gamma\in[0,1)$, reveals an unavoidable error in \ac{DMD} reconstruction. Sampling the solution with fixed step size $dt$, we get a one dimensional data matrix. Thus, the only valid \ac{DMD} is when the dimensionality reduction is maximal ($r=1$). In that case, the \ac{DMD} error (for details, see Appendix A, Eq. \eqref{eq:DMDERR}) converges to zero as $dt\to 0$. However, the reconstruction error (Eq. \eqref{eq:errorRec} is bounded from below.
A solution to this problem, as suggested in \cite{cohen2021modes}, is to formulate a time rescaled PDE by homogeneity normalization and to apply \ac{DMD} on that flow. However, this solution is not valid for zero homogeneous flows other than very simple cases, leading to another flaw in \ac{DMD}.

\subsection{Discontinuous Dynamical Systems}
The analytic solution of \eqref{motiEq:PDE} for $\gamma = 0$ is known when the Banach space is $\mathbb{R}^N$  . Applying the homogeneity normalization on zero homogeneous flow,  we find discontinuity in the dynamical modes. Thus, \ac{DMD} is not valid when the modes vanishes in finite time. The time rescale of \eqref{motiEq:PDE} and the relation of \ac{DMD} to zero-homogeneous decomposition is detailed in \cite{cohen2021Total}.

\subsection{Eigenvalue Multiplication}
\ac{DMD} is an exponential data fitting algorithm, \cite{askham2018variable}. Thus, \ac{DMD} can recover precisely the dynamics only when the typical decay profile of the system is exponential. However, even for the limited case of exponential decays, \ac{DMD} is not guaranteed to recover the dynamics. Let us consider a dynamic with a solution of the form
\begin{equation}
    u(t)=v\left(e^{\lambda_1 t}+e^{\lambda_2 t}\right).
\end{equation}
This solution cannot be reconstructed by a linear decomposition since the mode $v$ is associated with two eigenvalues ,$\lambda_1$ and $\lambda_2$.

The rest of this monograph attempts to propose a comprehensive solution to the aforementioned problems. We analyze the conditions for the existence of Koopman eigenfunctions and formulate the \ac{KMD} modes. Since \ac{DMD} is an approximation of  \ac{KMD}, if the \acp{KEF} do not exist the  approximation with \ac{DMD} is meaningless. After formulating the \ac{DMD} limitations we propose an alternative mode decomposition, which coincides with \ac{KMD} modes in a much broader setting.

\section{Koopman Eigenfunctions and Modes}\label{sec:disKoopman}
Koopman theory provides a linear representation to nonlinear dynamics by defining a new coordinate system. These coordinates are the measurements in the state-space termed as \emph{Koopman eigenfunctions}. Necessary and sufficient conditions for their existence are formulated here. Since the eigenfunctions are not unique, we define the \emph{Koopman family}, an infinite set of Koopman eigenfunctions. We also define a useful notion, referred to as the \emph{ancestors of a Koopman family}. This allows the reconstruction of the dynamical system, under certain conditions. Moreover, it allows to considerably enlarge the \acf{ROA} of the system. The above conclusions are consequences of the attributes of the dynamics, $P$, and its solution, $\bm{x}(t)$, discussed and analyzed below.

\subsection{Koopman Eigenfunctions}
\begin{figure}[htpb]
    \centering
    \includegraphics[width=0.85\textwidth,valign=c]{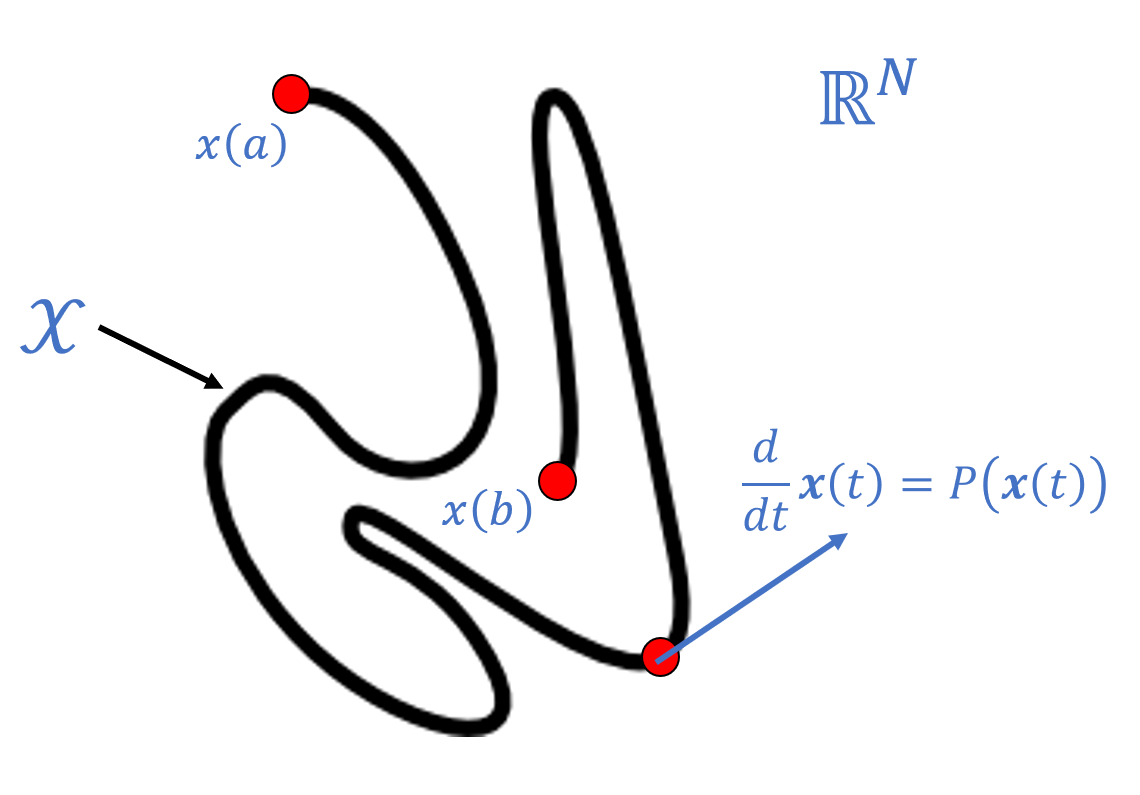}
    \caption[A curve]{The dynamic solution is represented as a curve\footnotemark}
    \label{fig:curve}
\end{figure}

We first set the necessary degree of smoothness of $P$,  required to develop the theory. This setting is highly non-restrictive and accommodates most useful linear and nonlinear dynamics, for both local and nonlocal settings. 
We refer to the operator $P$ in the dynamical system \eqref{eq:disDynamicalSystem}.
\begin{assumption}[Piecewise Continuous $P$]\label{assu:c0dynamics}
The operator $P:\mathbb{R}^N\to \mathbb{R}^N$ is in $C^0$ a.e. with zero Dirac measures.
\end{assumption}
This leads to the following Lemma.
\begin{lemma}[Continuous solution $\bm{x}(t)$]\label{lemma:continuousSolution}
If the operator $P$ in  Eq. \eqref{eq:disDynamicalSystem} admits  Assumption \ref{assu:c0dynamics} then the solution, $\bm{x}$, is in $C^1$ a.e.
\end{lemma}
\begin{proof}
The solution of the dynamics is
\begin{equation}
   \bm{x}(t)=\bm{x}(a)+\int_a^tP(\bm{x}(\tau))d\tau.
\end{equation}
The solution, $\bm{x}(t) \in C^1\, a.e.$  since $P\in C^0 \,a.e.$ and does not contain Dirac measures. 
\end{proof}
The solution, $\bm{x}(t)$, $t \in I\subset \mathbb{R}^+$, maps from the time range $I$ to $\mathbb{R}^N$. It can be interpreted as a parametric curve in $\mathbb{R}^N$, where its tangential velocity is $P(\bm{x})$. Let us denote the image of $\bm{x}(t)$ as $\mathcal{X}$. The image is the path in $\mathbb{R}^N$ where the system passes along the interval $I$. In Fig. \ref{fig:curve} an illustration of the solution of a dynamical system is shown. Using the Kinematics analogy, we can say the dynamics is a mass going from $\bm{x}(a)$ to $\bm{x}(b)$ with the instantaneous velocity, $P(\bm{x}(t))$, for every $t\in I$.
We note that Lemma \ref{lemma:continuousSolution} holds also if Assumption \ref{assu:c0dynamics} is limited to $\mathcal{X}$. 

\footnotetext{The curve image is taken from \url{https://mathinsight.org/definition/simple_curve} } 

A Koopman eigenfunction is a measurement of the solution $\bm{x}$ that admits Eq. \eqref{eq:KEFder} on the curve $\mathcal{X}$. As recently was stated in \cite{bollt2021geometric}, a Koopman eigenfunction can be formulated as an exponential function, where its argument is the inverse mapping from $\mathcal{X}$ to $I$ . The formal definition of the mapping is as follows.
\begin{definition}[Time state-space mapping]\label{def:invStaSpa}
    Let $\bm{x}(t)$ be the solution of the dynamical system \eqref{eq:disDynamicalSystem} where $t\in I$. Let  $\xi:\mathcal{X}\to I$ be a time state-space mapping from $\bm{x}$ to $t$, 
    \begin{equation}\label{eq:timeStaSpa}
        t= \xi(\bm{x}).
    \end{equation}
\end{definition}
This mapping is possible if the curve $\mathcal{X}$ is simple and open. Necessary conditions of a curve to be simple are discussed, for instance, in \cite{chuaqui2018general} and the references therein. 

\begin{lemma}[Differentiation of time state-space mapping ]\label{lem:timeDifferentiation} Let the conditions of Lemma \ref{lemma:continuousSolution} hold. If the time state-space mapping, $t=\xi(\bm{x})$, exists then it admits the following,
\begin{equation}\label{eq:coroTD}
    \nabla\xi(\bm{x})^T P(\bm{x})=1\quad a.e.\,\, \textrm{in}\,\, \mathcal{X}.
\end{equation}
\end{lemma}
\begin{proof}
The mapping  $\xi(\bm{x})$ is in $C^1$ a.e. in $\mathcal{X}$ since $\bm{x}(t)\in C^1 \, a.e.$ in $I$. The time derivative of the mapping is,
\begin{equation}
    1=\frac{d}{dt}t=\frac{d}{dt}\xi(\bm{x})=\nabla\xi(\bm{x})^T\frac{d\bm{x}}{dt}=\nabla\xi(\bm{x})^TP(\bm{x}).
\end{equation}
This expression is valid almost everywhere.
\end{proof}

We now turn to discuss necessary and sufficient conditions for the existence of a nontrivial Koopman eigenfunction (that is, a nonzero function which admits Eq. \eqref{eq:KEFder} with $\lambda\ne0$).

\begin{proposition}[Condition for the inexistence of a Koopman eigenfunction]\label{prop:equilibriumPoint}
If there is an equilibrium point in $I$ then a nontrivial Koopman eigenfunction does not exist.
\end{proposition}
\begin{proof}
Let $t_0\in I$ be an equilibrium point and $\varphi(\bm{x}(t))$ be a Koopman eigenfunction. Then, $\bm{x}(t)=const,\,\forall t\in[t_0,b]$. Therefore, Eq. \eqref{eq:KEFder} does not hold for nontrivial $\varphi$ for any $\lambda \ne 0$.
\end{proof}

\begin{remark}[Finite support time dynamics]\label{rem:inexist}
Let $P(\bm{x})$ define a dynamic for which the solution has a finite support in time. Namely, there is an extinction time point, $T_{ext}$, for which $P(\bm{x}(t))=\bm{0},\,\forall t\ge T_{ext}$. Then, if $T_{ext}\in I$, a Koopman operator $K_P^\tau$ has no eigenfunctions. We observe here that the time interval $I$ is crucial for the existence or inexistence of eigenfunctions.
\end{remark}

From a differential geometry perspective, as noted above, $\bm{x}(t)$ forms a curve where its tangential velocity is $P(\bm{x})$. The absence of an equilibrium point is equivalent to nonzero velocity. This type of parametric curves, where the velocity is always nonzero, is called \emph{regular}. The  Koopman eigenfunction does not exist for non-regular curves. 

\begin{lemma}[Koopman Eigenfunctions induced by a time state-space mapping]\label{lem:mapping2KE}
Let the conditions of Lemma \ref{lemma:continuousSolution} hold and $\bm{x}(t)$ be the solution of Eq. \eqref{eq:disDynamicalSystem}. If there exists a time state-space mapping, $t=\xi(\bm{x})$, then a Koopman eigenfunction exists a.e. in $I$.
\end{lemma}
\begin{proof} The mapping, $t=\xi(\bm{x})$, is in $C^1$ a.e. in $\mathcal{X}$ since $\bm{x}(t)$ is in $C^1$ a.e. in $I$. Given that mapping, we define the following function,
\begin{equation}\label{eq:proofLemTM}
    \varphi(\bm{x})=e^{\alpha \xi(\bm{x})+\beta}.
\end{equation}
This function is in $C^1$ a.e. in $\mathcal{X}$. The time derivative of this function is, 
\begin{equation}
\begin{split}\label{eq:timeStateDet}
    \frac{d}{dt}\varphi(\bm{x}(t))&=\frac{d}{d\xi} e^{\alpha\cdot \xi(\bm{x})+\beta}\nabla\xi(\bm{x})^T \frac{d}{dt}\bm{x}=\alpha \varphi(\bm{x}(t))\nabla\xi(\bm{x})^T P(\bm{x}(t)).
\end{split}
\end{equation}
According to Lemma \ref{lem:timeDifferentiation}, $\nabla\xi(\bm{x})^T P(\bm{x}(t))=1$ a.e.. Thus, the function in Eq. \eqref{eq:proofLemTM} admits Eq. \eqref{eq:KEFder} for any value of $\beta$, where the corresponding eigenvalue is $\lambda = \alpha$. 
\end{proof}

\begin{theorem}[Sufficient condition for the existence of a Koopman eigenfunction] \label{theo:sufficientConditionKEF} Let the conditions of Lemma \ref{lemma:continuousSolution} hold and one of the entries of the vector $P(\bm{x}(t))$ 
is either positive or negative $\forall t \in I$. Then, Koopman eigenfunctions exist a.e. in the time interval $I$. 
\end{theorem}
\begin{proof}
If one of the entries in $P(x(t))$ is either positive or negative for all $t\in I$ then this entry is monotone and therefore injective. Then, the curve $\mathcal{X}$ is simple and open (see \cite{courant2012introduction} pages 45, 177 and 207). Therefore, the time state-space mapping, $\xi(\bm{x})$, exists.
Following Lemma \ref{lem:mapping2KE}, Koopman eigenfunctions can be expressed by \eqref{eq:proofLemTM}.
\end{proof}

The simple example below illustrates the connections between the equilibrium point, finite time dynamics and time state-space mapping.
\begin{idoexample}[Finite time support]\label{example:FiniteTimeSupport}
Let us consider the following dynamics,
\begin{equation}
    \frac{d}{dt}x=-2x^{\frac{1}{2}},\quad x(0)=1.
\end{equation}
The solution is
\begin{equation}
    x(t) = \begin{cases}(1-t)^2&t\in[0,1]\\
    0&t>1
    \end{cases}.
\end{equation}
For $I=[0,1]$, the time state-space mapping is,
\begin{equation}\label{eq:exampleFDTimeMapping}
    t(x) = 1-\sqrt{x}, \quad I=[0,1],
\end{equation}
and using \eqref{eq:proofLemTM} with $\alpha=1$, $\beta=0$, we can express a Koopman eigenfunction by, 
\begin{equation}
    \varphi(x)=e^{1-\sqrt{x}}.
\end{equation}
Now, let us repeat this example with a different time interval. Let $I=[0,1.5]$, containing the extinction time $T_{ext}=1$. Note that, first, the time mapping, Eq. \eqref{eq:exampleFDTimeMapping}, does not hold in the entire interval, and the eigenfunction $\varphi$ does not admit $\frac{d}{dt}\varphi(x)=\varphi(x)$ since $\varphi(x)$ is a nonzero constant for $t\in[1,1.5]$. 
\end{idoexample}

\subsubsection{Extended \texorpdfstring{\ac{DMD}}{TEXT} induced from time state-space mapping}
One of the methods to increase the accuracy of the classic \ac{DMD} is by enriching the state-space vector with nonlinear measurements of the coordinates $\bm{x}$, see \cite{williams2015data}. 
It is shown that this approach indeed improves accuracy, however -  the theoretical justification is vague. In addition - the enriching method may become somewhat heuristic.
We can interpret this approach as the Taylor expansion of Koopman eigenfunctions. 
This provides both justification and a clear method for supplying additional measurements.
Let us expand the Koopman eigenfunction, $\varphi(\bm{x})=e^{\xi(\bm{x})}$, by a Taylor series,
\begin{equation*}
    \varphi(\bm{x})=e^{\xi(\bm{x})}=\sum_{j=0}^\infty \frac{\xi(\bm{x})^j}{j!}.
\end{equation*}
We can approximate this expression by taking only finite number of elements from this sum,
\begin{equation*}
    e^{\xi(\bm{x})}\approx\sum_{j=0}^M \frac{\xi(\bm{x})^j}{j!}.
\end{equation*}
Thus, Eq. \eqref{eq:KEFder} can be approximated as
\begin{equation}\label{eq:TaylorKEF}
    \frac{d}{dt}\sum_{j=0}^M \frac{\xi(\bm{x})^j}{j!}\approx\sum_{j=0}^M \frac{\xi(\bm{x})^j}{j!}.
\end{equation}
In matrix notation, this approximation can be reformulated as,
\begin{equation}\label{eq:EDMD}
\frac{d}{dt}
\begin{bmatrix}
1\\\xi(\bm{x})\\ \vdots\\\frac{\xi(\bm{x})^M}{M!}
\end{bmatrix}\approx A \begin{bmatrix}
1\\\xi(\bm{x})\\ \vdots\\\frac{\xi(\bm{x})^M}{M!},
\end{bmatrix}
\end{equation}
where any matrix $A$ with a left-eigenvector $\begin{bmatrix}
1&\cdots&1
\end{bmatrix}$ can be an optional solution to Eq. \eqref{eq:EDMD} for which Eq. \eqref{eq:TaylorKEF} holds. In addition, taking $M$ to infinity, $A$ gets the form
\begin{equation*}
    [A]_{i,j}=\begin{cases}
        1& i=j+1\\
        0& else
    \end{cases}
\end{equation*}
where $[A]_{i,j}$ is the ${i,j}$th entry of $A$.

\subsection{Koopman Family}
The \ac{KEF} is of the form $\varphi(t)=e^{\alpha t+ \beta}$, Eq. \eqref{eq:proofLemTM}. This form of solution is unique following a standard existence and uniqueness theorem of ODE's. The exponential parameters, $\alpha$ and $\beta$, are dictated by the eigenvalue and the initial condition. Without these restrictions, there are infinite \acp{KEF} for any dynamical system.

From a different angle, viewing the state-space $\bm{x}$ as a curve gives a compelling interpretation of the multiplicity of Koopman eigenfunctions. A curve can be reparameterized in different manners. Changing the parameters, $\alpha$ and $\beta$, amounts to reparameterization by translation and scaling. This insight leads us to the following lemma, which extends the identities presented in \cite{bollt2021geometric}. We show that any mathematical manipulation on a \acs{KEF} which maintains the form of Eq. \eqref{eq:proofLemTM} generates a new \ac{KEF}.

\begin{lemma}[Multiplicities of Koopman eigenfunctions] \label{lemma:KoomanFamily} If $\varphi_1, \varphi_2$ are Koopman eigenfunctions with the corresponding eigenvalues $\lambda_1,\lambda_2$ then:
\begin{enumerate}[leftmargin=10mm]
    \item The function $a\cdot \varphi_1, \, a\in\mathbb{R}$ ($a\ne 0$) is an eigenfunction with eigenvalue $\lambda_1$. 
    \item The function $(\varphi_1)^\alpha, \, \alpha\in\mathbb{C}$ ($\alpha\ne 0$) is an eigenfunction with eigenvalue $\alpha\lambda_1$. 
    \item For any $n,m\in \mathbb{R}$  the function $\left(\varphi_1\right)^n\left(\varphi_2\right)^m$ is  an eigenfunction with eigenvalue $n\lambda_1+m\lambda_2$.
    \item The function $\left(\varphi_1\right)^\frac{\lambda}{\lambda_1}+ \left(\varphi_2\right)^\frac{\lambda}{\lambda_2}$ is  an eigenfunction with eigenvalue $\lambda$.
\end{enumerate}
\end{lemma}
\begin{proof}
This can be shown by,
\begin{enumerate}[leftmargin=10mm]
    \item Using the linearity of the Koopman operator.
    \item Writing the time derivative of $(\varphi_1)^\alpha$ explicitly we get, 
    \begin{equation}
        \begin{split}
            \frac{d}{dt}\left[\varphi_1^\alpha\right]&=\alpha(\varphi_1)^{\alpha-1}\lambda_1\varphi_1=\alpha \lambda_1 \varphi_1^\alpha.
        \end{split}
    \end{equation}
    \item Similarly, 
    \begin{equation}
        \begin{split}
            \frac{d}{dt}\left[(\varphi_1)^n(\varphi_2)^m\right]&=(\varphi_2)^m n(\varphi_1)^{n-1}\lambda_1\varphi_1 +(\varphi_1)^n m(\varphi_2)^{m-1}\lambda_2 \varphi_2\\
            &=(n\lambda_1+m\lambda_2)(\varphi_1)^n(\varphi_2)^m.
        \end{split}
    \end{equation}
    
    \item Finally, 
    \begin{equation}
        \begin{split}
            \frac{d}{dt}\left[\left(\varphi_1\right)^\frac{\lambda}{\lambda_1}+ \left(\varphi_2\right)^\frac{\lambda}{\lambda_2}\right]&=\frac{\lambda}{\lambda_1}\left(\varphi_1\right)^{\frac{\lambda}{\lambda_1}-1}\lambda_1 \varphi_1+ \frac{\lambda}{\lambda_2}\left(\varphi_2\right)^{\frac{\lambda}{\lambda_2}-1}\lambda_2 \varphi_2\\
            &=\lambda\left[\left(\varphi_1\right)^\frac{\lambda}{\lambda_1}+ \left(\varphi_2\right)^\frac{\lambda}{\lambda_2}\right].
        \end{split}
    \end{equation}
\end{enumerate}
\end{proof}
\paragraph{Discussion} The multiplicities presented in Lemma \ref{lemma:KoomanFamily} are translation and scaling of the time variable. Case 1 in this Lemma is a translation of the time axis and the rest of the cases are scaling. To distinguish between eigenfunctions which are generated from other eigenfunctions and ``new'' independent ones, we introduce the concepts of \emph{Koopman family} and its \emph{ancestors}.

\begin{definition}[Koopman family]\label{def:koopmanFamily}
    Let $\{\varphi_i\}_{i=1}^n$ be a finite set of \acp{KEF}. Let $k_P(\{\varphi_i\}_{i=1}^n)$ be the infinity uncountable set of \acp{KEF} generated by the finite set, recursively, according to the four options stated in Lemma \ref{lemma:KoomanFamily}. Let us define $k_P^m(\{\varphi_i\}_{i=1}^n)=k_P(k_P^{m-1}(\{\varphi_i\}_{i}^n))$. We term $\mathcal{K}_P(\{\varphi_i\}_{i}^n)= \cup_{m=1}^\infty k_P^m(\{\varphi_i\}_{i=1}^n)$ as the Koopman family of $\{\varphi_i\}_{i=1}^n$.
\end{definition}

\begin{definition}[Ancestors of a Koopman family]\label{def:koopmanFamilyAncestors}
    Let $\{\varphi_i^*\}_{i=1}^m$ be a finite set of \acp{KEF}. This set is an ancestor set of the Koopman family $\mathcal{K}_P(\{\varphi_i\}_{i=1}^n)$ if the following conditions hold:
    \begin{enumerate}[leftmargin=10mm]
        \item $\varphi\in \mathcal{K}_P(\{\varphi_i\}_{i}^n)\Longleftrightarrow \varphi\in \mathcal{K}_P(\{\varphi_i^*\}_{i}^m)$.
        \item $\varphi_j^* \notin \mathcal{K}_P(\{\varphi_i^*\}_{i=1,i\ne j}^m)$ for any $j=1,2,\cdots,m$.
    \end{enumerate}
\end{definition}
Note that the subscript $_P$ is for the dynamical system.

\subsubsection{Koopman Eigenfunction Vector}
A vector of Koopman eigenfunctions is denoted by,
\begin{equation}\label{eq:KEFVectorDef}
    \bm{\varphi}(\bm{x})=\begin{bmatrix}
    \varphi_1(\bm{x})&\cdots&\varphi_L(\bm{x})
    \end{bmatrix}^T,
\end{equation}
where $L$ can be finite or infinite. The Jacobian matrix of this vector is,
\begin{equation}\label{eq:KEFVectorDefJacob}
    \frac{\partial}{\partial \bm{x}} \bm{\varphi}(\bm{x})=\begin{bmatrix}\nabla \varphi_1(\bm{x})^T\\
    \vdots\\
    \nabla\varphi_L(\bm{x})^T\end{bmatrix}=\mathcal{J}(\bm{\varphi}).
\end{equation}

\begin{theorem}[Linear dynamic in Koopman family]\label{theo:lini} Let the conditions of Theorem \ref{theo:sufficientConditionKEF} hold. The dynamical system $P$  can be represented as a linear one with a vector of Koopman eigenfunctions, where the time derivative of this vector is,
\begin{equation}
    \frac{d}{dt}\bm{\varphi}(\bm{x})=\mathcal{J}(\bm{\varphi})P(\bm{x})=\Lambda\bm{\varphi}(\bm{x}),\quad a.e.
\end{equation}
where $\Lambda$ is a diagonal matrix with the corresponding eigenvalues.
\end{theorem}
\begin{proof} We would like to prove first the existence of a $L$ dimensional \ac{KEF}. From Theorem \ref{theo:sufficientConditionKEF} there exists a \ac{KEF}. From Lemma \ref{lemma:KoomanFamily} if there exists a \ac{KEF}, there are infinite set of \acp{KEF}, therefore, at least $L$ eigenfunctions. According to the definition of the Koopman eigenfunction, the time derivative is,
\begin{equation}
    \frac{d}{dt} \bm{\varphi}(\bm{x})=\begin{bmatrix}\frac{d}{dt}\varphi_1(\bm{x}),\cdots,\frac{d}{dt}\varphi_L(\bm{x})\end{bmatrix}^T=\begin{bmatrix}\lambda_1\varphi_1(\bm{x}),\cdots,\lambda_L\varphi_L(\bm{x})\end{bmatrix}^T=\Lambda \bm{\varphi}(\bm{x}).
\end{equation}
On the other hand, applying the chain rule we get,
\begin{equation}
    \frac{d}{dt} \bm{\varphi}(\bm{x})=\begin{bmatrix}\nabla \varphi_1(\bm{x})^T\\
    \vdots\\
    \nabla\varphi_L(\bm{x})^T\end{bmatrix}\frac{d}{dt}\bm{x}(t)=\mathcal{J}(\bm{\varphi})P(\bm{x}).
\end{equation}
Note that $\varphi_i(\bm{x})\in C^1,\,(a.e.)$, so the expressions above are valid only almost everywhere.
\end{proof}

\subsubsection{Reconstructing the dynamics}
The ability to reconstruct the dynamics is based on the relations between the vectors ${\bm \varphi}$ and ${\bm x} $. In classical control theory, this is referred to as \emph{observability}. Here, we suggest to examine the notion of observability by computing the rank of the Jacobian matrix $\mathcal{J}(\bm{\varphi})$,  Eq. \eqref{eq:KEFVectorDefJacob}. The rows of this matrix are the gradients of the \acp{KEF}. The following lemma shows that the gradient of a member of a  Koopman family originates with its ancestors.

\begin{lemma}[KEF gradients of a family]\label{lemma:KEFGradients}
Let $\mathcal{K}_P(\{\varphi_i^*\}_{i=1}^m)$ be the Koopman family of an ancestor set, $\{\varphi_i^*\}_{i=1}^m$. Let $\varphi$ be a \ac{KEF} in $\mathcal{K}_P(\{\varphi_i^*\}_{i=1}^m)$. Then, the gradient of $\varphi$, $\nabla \varphi$, is a linear combination of the gradients of the ancestor set for any $t\in I$.
\end{lemma}
\begin{proof}
Let $\mathcal{KG}$ be the linear span, defined by
\begin{equation}
    \mathcal{KG}=span\left(\left\{\nabla \varphi_i^*\right\}_{i=1}^m\right)=\left\{\sum_{i=1}^m a_i\nabla \varphi_i^*,\,\, \forall a_i\in \mathbb{C}\right\}.
\end{equation}
Let  $\varphi$ be in  $\mathcal{K}_P(\{\varphi_i^*\}_{i=1}^m)$. According to Definition \ref{def:koopmanFamily}, there exist recursive steps leading from the ancestors $\{\varphi_i^*\}_{i=1}^m$ to $\varphi$. Now, by induction we show that $\nabla\varphi\in\mathcal{KG}$. Let us assume that from the ancestors to $\varphi$ there is one step. Namely, $\varphi$ is generated using  $\varphi^*_i,\varphi^*_j$, according to the four cases of Lemma \ref{lemma:KoomanFamily}. The gradient of $\varphi$ is a linear combination of the gradients of $\varphi^*_i$ and $\varphi^*_j$. For cases 1 ,2 and 4, the linearity is straightforward. For case 3, $\varphi=(\varphi^*_i)^n (\varphi^*_j)^l$, we have,
\begin{equation}
\begin{split}
    \nabla{\varphi}&=\nabla{\left(\varphi^*_i\right)^n \left(\varphi^*_j\right)^l}=n\left(\varphi^*_i\right)^{n-1}\left(\varphi^*_j\right)^l \nabla{\left(\varphi^*_i\right) }+l\left(\varphi^*_i\right)^n \left(\varphi^*_j\right)^{l-1}\nabla{\left(\varphi^*_j\right)}\\
    &=\begin{bmatrix}
    \nabla{\left(\varphi^*_i\right) }&\nabla{\left(\varphi^*_j\right)}
    \end{bmatrix}\begin{bmatrix}n\left(\varphi^*_i\right)^{n-1}(\varphi^*_j)^l\\
    l\left(\varphi^*_i\right)^n (\varphi^*_j)^{l-1}
    \end{bmatrix}.
    \end{split}
\end{equation}
For any $t\in I$ the vector $\begin{bmatrix}
n\left(\varphi^*_i\right)^{n-1}(\varphi^*_j)^l&
    l\left(\varphi^*_i\right)^n (\varphi^*_j)^{l-1}
\end{bmatrix}^T$ is constant. Therefore, the gradient of $\varphi$ is in $\mathcal{KG}$.
Now we assume there exist $k$ steps from the ancestors to $\varphi$. Let $\varphi_1$ and $\varphi_2$ be generated by $k-1$ steps. The induction assumption holds, meaning, their gradients are in $\mathcal{KG}$. Now, there is one step from $\varphi_1$ and $\varphi_2$ to $\varphi$. As shown,  $\nabla\varphi$ is a linear combination of the gradients of its generators, $\nabla \varphi_1, \nabla \varphi_2$. But these vectors belong to $\mathcal{KG}$ by the induction assumption. Therefore, $\nabla \varphi \in \mathcal{KG}$. 
\end{proof}

The multiplicity of Koopman eigenfunctions results from either arithmetical manipulations (Def. \ref{def:koopmanFamily}) or the existence of several time state-space mappings (Def. \ref{def:invStaSpa}). The main difference is the rank of the Jacobian, $\mathcal{J}(\bm{\varphi})$. Given a vector of \acp{KEF}, adding another Koopman eigenfunction from the Koopman family of the \acp{KEF} in the vector -- does not increase the rank of the Jacobian. However, adding a Koopman eigenfunction from another time state-space mapping does. The Jacobian matrix rank is related to system controllability and observability (see for example \cite{brunton2019data,evangelisti2011controllability}). In the following, we formulate the connections between the rank of  the Jacobian matrix, the size of the ancestor set, and time state-space mappings.

\begin{definition}[Full observability in the context of Koopman theory] Consider the dynamical system Eq. \eqref{eq:disDynamicalSystem} where $\bm{x}\in \mathbb{R}^N$. The system is fully observable if the state-space can be revealed from the \acp{KEF} and the initial condition.
\end{definition}
\begin{proposition}[Sufficient conditioin for full observability]\label{prop:FullObservability} Consider the dynamical system \eqref{eq:disDynamicalSystem} where $\bm{x}\in \mathbb{R}^N$. Let us denote the Koopman family of all Koopman eigenfunctions of the dynamics as $\mathcal{K}_P$. An ancestor set of $\mathcal{K}_P$ is denoted as $\{\varphi_i^*\}_{i=1}^n$. The system is fully observable if $N\le n$.
\end{proposition}
\begin{proof} According to Lemma \ref{theo:lini}, for any vector of \acp{KEF} the following equation holds,
\begin{equation}
    \mathcal{J}(\bm{\varphi})P(\bm{x})=\Lambda\bm{\varphi}(\bm{x}),\quad (a.e.).
\end{equation}
Let us choose a vector of ancestors, i.e.
\begin{equation}
    \bm{\varphi^*}(\bm{x})=\begin{bmatrix}
    \varphi_1^*(\bm{x})&\cdots&\varphi_n^*(\bm{x})
    \end{bmatrix}^T.
\end{equation}
According to Lemma \ref{lemma:KEFGradients} the rank of the Jacobian matrix is full and equal to $N$. Since the matrix $\mathcal{J}(\bm{\varphi})^T\mathcal{J}(\bm{\varphi})$ is invertible, the dynamics, $P$, can be revealed according to the following relation,
\begin{equation}\label{eq:systemRecovery1}
    P(\bm{x})=\left(\mathcal{J}(\bm{\varphi})^T\mathcal{J}(\bm{\varphi})\right)^{-1}\mathcal{J}(\bm{\varphi})^T\Lambda\bm{\varphi}(\bm{x}),\quad (a.e.).
\end{equation}
That is, we use the Moore-Penrose inverse. The state-space can now be calculated as,
\begin{equation}
    \bm{x}(t)=\bm{x}_0+\int_a^t\left(\mathcal{J}(\bm{\varphi}(\tau))^T\mathcal{J}(\bm{\varphi}(\tau))\right)^{-1}\mathcal{J}(\bm{\varphi}(\tau))^T\Lambda\bm{\varphi}(\tau)d\tau.
\end{equation}
\end{proof}

\begin{corollary}[Full observability for a monotone  dynamics] \label{corollary:fullyObservableSystem}
If each entry in $P$ is either positive or negative for any $t$ in  $I$ then each entry of the state-space is monotone (and injective). We can formulate $N$ different time state-space mappings from $\mathcal{X}$ to $I$ (Theorem \ref{theo:sufficientConditionKEF}). These mappings induce $N$ different \acp{KEF} and according to Proposition \ref{prop:FullObservability} the system is fully observable. 

\end{corollary}
\begin{remark}[Sufficient condition for dynamic reconstruction]\label{remark:sufficientDynamicRecon}
If each of the entries of $P$ is either positive or negative for all $t$ in $I$ then the dynamics can be reconstructed as
\begin{equation}\label{eq:dynamicReconstruction}
    P(\bm{x})=\mathcal{J}^{-1}(\bm{\varphi})\Lambda \bm{\varphi}(\bm{x}).
\end{equation}
According to Corollary \ref{corollary:fullyObservableSystem}, if each of the entries of $P$ is either positive or negative in $I$ then the Jacobian matrix is $N\times N$ and is full rank, therefore -- invertible. Using Theorem \ref{theo:lini} we reach Eq. \eqref{eq:dynamicReconstruction}.
\end{remark}

\begin{remark}[Global controllability]\label{remark:GlobalControllability} Reconstructing the dynamical system enables us to enlarge the \acf{ROA}, Eq. \eqref{eq:ROA}. Given the nonlinear dynamics,
\begin{equation}
    \frac{d}{dt}\bm{x}(t)=P(\bm{x}(t))+\bm{u},
\end{equation}
we can cancel the nonlinearity with the ancestors of a Koopman family $\bm{\varphi}(\bm{x})$ if the dynamics is fully observable. In order to reach a stable system for any point $\bm{x}$ we define the following input $\bm{u}$, 
\begin{equation}
    \bm{u}=\mathcal{J}^{-1}(\bm{\varphi})\Lambda \bm{\varphi}(\bm{x})+\bm{w}, 
\end{equation}
where the first element cancels the nonlinearity of the system (Remark \ref{remark:sufficientDynamicRecon}) and the second term brings the system to any  desired point in $\mathbb{R}^N$. Note that we assume there are no singular points in  $P$.
\end{remark}

\begin{remark}[Reconstructing the dynamic - limitations] The recovery of the system, as described by Eq. \eqref{eq:systemRecovery1} and \eqref{eq:dynamicReconstruction}, is valid for a given initial condition $x_0$. In order to obtain a full recovery of the system, the properties of the initial condition should be taken into account. These subject exceeds the frame of this work and requires further research.
\end{remark}

\subsubsection{Reconstruction conservation laws}
Dynamic reconstruction and conservation laws (such as energy, momentum etc.) are perhaps the most crucial tasks  in dynamical system analysis, in general, and controlling systems in particular. Data driven algorithms to reveal the dynamic (governing laws) and the conservation laws, based on the Koopman operator theory have been studied by \cite{rudy2017data,brunton2016discovering,schmidt2009distilling,kaiser2018discovering,langley1981bacon}. The common approach argues that the conservation laws are related to the null part of the Koopman spectrum. Namely, the Koopman eigenfunctions related to eigenvalue zero are or may be formulations of the conservation laws. In a similar manner, the dynamic can be reconstructed based on the nontrivial Koopman eigenfunctions. 

We propose an alternative view. As discussed above, the relevant \acp{KEF} to system reconstruction are indeed not in the null part of the Koopman spectrum. However, from our perspective, after recovering the dynamical system (the governing laws) via \acp{KEF} the conservation laws naturally emerge from these nontrivial \acp{KEF}.

Let $\varphi(x)$ be an eigenfunction, namely Eq. \eqref{eq:KEFder} holds for some $\lambda$. We consider the measurement $\ln\left(\varphi(x)\right)$. By using Eq. \eqref{eq:KEFdef}, we can express the time derivative of this measurement as,
\begin{equation*}
    \frac{d}{dt}\ln\left(\varphi(\bm{x})\right)=\frac{\lambda\varphi(\bm{x})}{\varphi(\bm{x})}=\lambda.
\end{equation*}
On the other hand, with the chain rule we get,
\begin{equation*}
    \frac{d}{dt}\ln\left(\varphi(\bm{x})\right)=\frac{\frac{d}{dt}\varphi(\bm{x})}{\varphi(\bm{x})}=\frac{\nabla \varphi(\bm{x})^T\frac{d}{dt}\bm{x}}{\varphi(\bm{x})}=\frac{\nabla \varphi(\bm{x})^TP(\bm{x})}{\varphi(\bm{x})}.
\end{equation*}
Then, for any Koopman eigenfunction (an intrinsic coordinate) we can formulate the following conservation law,
\begin{equation}\label{eq:rawConLaw}
    \frac{\nabla \varphi(\bm{x})^TP(\bm{x})}{\varphi(\bm{x})}=\lambda.
\end{equation}
This result coincides with Eq. \eqref{eq:coroTD}, by using Lemma \ref{lem:timeDifferentiation} and formulating a \ac{KEF} by a time state-space mapping, Eq. \eqref{eq:proofLemTM}, with $\alpha =\lambda$. In other words, when a time state-space mapping $\xi$ exists, an alternative formulation to the conservation law of Eq. \eqref{eq:rawConLaw} is,
\begin{equation}
    1=\frac{d}{dt}t=\frac{d}{dt}\xi(\bm{x})=\nabla\xi(\bm{x})^T\frac{d\bm{x}}{dt}=\nabla\xi(\bm{x})^TP(\bm{x}).
\end{equation}
We illustrate this with the following two examples.

\begin{idoexample}[Free Fall]\label{exm:freeFall} Let $x_1(t)$ and $x_2(t)$ be the height and the velocity of a mass in a free fall, respectively. The dynamical system is,
\begin{equation*}
\frac{d}{dt}\begin{bmatrix}
x_1\\x_2
\end{bmatrix}=\begin{bmatrix}
x_2\\-g
\end{bmatrix},
\end{equation*}
where the initial condition is $\begin{bmatrix}
h&0
\end{bmatrix}^T$. The solution is 
\begin{equation*}
\begin{bmatrix}
x_1\\x_2
\end{bmatrix}=\begin{bmatrix}
h-\frac{1}{2}gt^2\\-gt
\end{bmatrix}.
\end{equation*}
The time state-space mappings are,
\begin{equation*}
    \begin{bmatrix}
    t_1\\t_2
    \end{bmatrix}=\begin{bmatrix}
    \sqrt[]{\frac{2(h-x_1)}{g}}\\
    -\frac{x_2}{g}
    \end{bmatrix}.
\end{equation*}
The induced conservation laws, using Eq. \eqref{eq:rawConLaw} are as follows.\\
{\bf{Conservation law \#1}}
\begin{equation*}
    1=\frac{d}{dx_1}t_1(x_1)\cdot \frac{d}{dt}x_1=-\frac{1}{2\sqrt[]{\frac{2(h-x_1)}{g}}}\frac{2}{g}\cdot x_2
\end{equation*}
One can reformulate this to the energy conservation law, 
\begin{equation*}
    gx_1+\frac{1}{2}x_2^2=hg.
\end{equation*}

{\bf{Conservation law \#2}}
\begin{equation*}
    \frac{d}{dx_2}t_2(x_2)\cdot \frac{d}{dt}x_2=-\frac{1}{g}\cdot (-g)=1
\end{equation*}
The conservation law \#2 is due to the constant acceleration, $g$.
\end{idoexample}

\begin{idoexample}[Pure rolling down an incline]  On an inclined plane with a slope of angle $\alpha$, a solid cylinder with mass $m$, radius $R$, and rotational inertia $I_{cm}$
is released from rest. The location along the plane is denoted by $x_1$ and its velocity by $x_2$. The dynamical system is,
\begin{equation}
    \frac{d}{dt}\begin{bmatrix}
    x_1\\x_2
    \end{bmatrix}=\begin{bmatrix}
    x_2\\
    \frac{g sin \alpha}{1+\frac{I_{cm}}{mR^2}}
    \end{bmatrix}
\end{equation}
with the initial condition $\bm{x}=[0,0]^T$. The solution is,
\begin{equation}
\begin{split}
    x_1(t)&=\frac{1}{2}\frac{g sin \alpha}{1+\frac{I_{cm}}{mR^2}}t^2,\\
    x_2(t) &= \frac{g sin \alpha}{1+\frac{I_{cm}}{mR^2}}t.
\end{split}
\end{equation}
The time mappings are 
\begin{equation}
\begin{split}
    t_1(x_1)&=\sqrt[]{2\frac{1+\frac{I_{cm}}{mR^2}}{g sin \alpha}x_1},\\
    t_2(x_2) &= \frac{1+\frac{I_{cm}}{mR^2}}{g sin \alpha}x_2.
\end{split}
\end{equation}
{\bf{Conservation law \#1}}
\begin{equation}
    1=\frac{d}{dx_1}t_1(x_1)\frac{d x_1}{dt}=\sqrt[]{2\frac{1+\frac{I_{cm}}{mR^2}}{g sin \alpha}}\frac{1}{2\sqrt{x_1}}x_2
\end{equation}
We can reformulate it as,
\begin{equation}
    \underbrace{\frac{1}{2}mx_2^2}_{E_K}+\underbrace{\frac{1}{2}I_{cm}\left(\frac{x_2}{R}\right)^2}_{E_R}\underbrace{-mgx_1sin\alpha}_{E_P} =0
\end{equation}
getting,  as expected, that the sum of the energies (Kinetic, Rotational, and Potential) is zero.

{\bf{Conservation law \#2}}
In the same manner as in Example \ref{exm:freeFall}, conservation law \#2 is a result of constant acceleration.
\end{idoexample}

\subsection{Koopman Mode Decomposition}
The Koopman mode decomposition leverages this infinite family to reconstruct the observations from the Koopman eigenfunctions (\cite{mezic2005spectral}). The reconstruction is a linear combination of Koopman eigenfunctions. For instance, the $i$th entry of $\bm{x}$ is assumed to be reconstructed as (\cite{brunton2021modern}),
\begin{equation}
    x_i(t)=\sum_{j=1}^{\infty}v_{i,j}\varphi_j(\bm{x}(t)),
\end{equation}
where $v_{i,j}$ is a scalar. Then, the state-space can be written as,
\begin{equation}\label{eq:reconKD}
    \bm{x}(t)=\sum_{j=1}^{\infty}\bm{v}_j\varphi_j(\bm{x}(t)),
\end{equation}
where $\bm{v}_j$ is an $N$ dimensional vector whose entries are the coefficients of the $j$th Koopman eigenfunction, namely $\bm{v}_j =\begin{bmatrix}
v_{1,j} & \cdots& v_{N,j}
\end{bmatrix}^T$. Substituting the solution of $\varphi(\bm{x})$, Eq. \eqref{eq:KEF}, we get,
\begin{equation}
    \bm{x}(t)=\sum_{j=1}^{\infty}\bm{v}_j\varphi_j(\bm{x}(a))e^{\lambda_j t}.
\end{equation}
The infinite triplet $\{\bm{v}_j,\varphi_j,\lambda_j\}_{j=1}^\infty$ is the Koopman mode decomposition, where $\{\bm{v}_j\}_{j=1}^\infty$ are the Koopman modes, $\{\varphi_j\}_{j=1}^\infty$ are the \acp{KEF}, and $\{\lambda_j\}_{j=1}^\infty$ are the Koopman eigenvalues. Note that the maximal index argument in the sum of Eq. \eqref{eq:reconKD} is not necessarily infinity. For example, it is enough to have one mode to reconstruct the linear dynamics initiated with one of its eigenvectors. In matrix notations, let $V$ be a matrix whose column vectors are the corresponding Koopman modes. The state-space can be expressed as,
\begin{equation}\label{eq:matModeDecomp}
    \bm{x}(t)=V\bm{\varphi}(\bm{x}(t)).
\end{equation}
Thus, the dynamical system has a linear representation with the measurements $\{\varphi_j(\bm{x})\}_{j=1}^\infty$, \cite{kaiser2021data}.
\begin{idoexample}[\ac{KMD} of Zero Homogeneous Dynamics] \label{exm:ZeroHomoDynamics}
Let us consider the following dynamical system
\begin{equation}
    \frac{d}{dt}\bm{x}= P(\bm{x}), \quad \bm{x}(t=0)=\bm{v},\, I=[0,-1/\lambda)
\end{equation}
where $P$ is a zero homogeneous operator (admitting $P(a\cdot \bm{x})=sign(a)P(\bm{x}),\,\forall a\in \mathbb{R}$), $\bm{v}$ and $\lambda$ 
are a nonlinear eigenvector and  the corresponding eigenvalue of $P$, respectively, i.e. they admit the nonlinear eigenvalue problem $P(\bm{v})=\lambda \bm{v}$. 
We assume a stable system, where $\lambda < 0$.
More background on such problems is presented in \cite{gilboa2018nonlinear}. Then, the solution of the ODE is,
\begin{equation}\label{eq:exmKMDZeroHomo}
    \bm{x}(t)=\bm{v}\left(1+\lambda t\right),\quad t\in I.
\end{equation}
A \ac{KEF}  can be formulated by the time state-space mapping as,
\begin{equation}
    \varphi(t)=e^t=e^{\frac{\frac{\inp{\bm{x}}{\bm{v}}}{\norm{\bm{v}}^2}-1}{\lambda}}.
\end{equation}
We would like now to express the solution \eqref{eq:exmKMDZeroHomo} with Koopman eigenfunctions. To express the function $t$ we have to apply the natural logarithm, $\ln$, on the Koopman eigenfunction. With  Taylor  series one can express it as,
\begin{equation}
    t=\ln (\varphi(\bm{x}))=\sum_{n=1}^\infty (-1)^{n+1}\frac{(\varphi(\bm{x})-1)^n}{n}.
\end{equation}
Then, the solution of \eqref{eq:exmKMDZeroHomo} can be written as,
\begin{equation}\label{eq:exmSolutionAsSumKEF}
    \bm{x}(t)=\bm{v}\left(1+\lambda \sum_{n=1}^\infty (-1)^{n+1}\frac{(\varphi(\bm{x})-1)^n}{n}\right).
\end{equation}
By expanding the terms $(\varphi - 1 )^n$ we get an infinite polynomial with respect to the \ac{KEF} $\varphi$. \ac{KMD} emerges naturally.
\end{idoexample}
\paragraph{Discussion} According to this example, since there is only one mode and its decay profile is not exponential, there can be many \acp{KEF} for one Koopman mode. The  multiplicity of eigenvalues for one mode is related to the limitations of \ac{DMD}. Since  \ac{DMD} recovers only linear dynamics it cannot handle well one eigenvector with multiple eigenvalues. 

We can now formulate the relation between Koopman modes and the dynamical system.
\begin{proposition}[The Jacobian and Koopman modes]\label{prop:JacKM}
Let $\bm{\varphi}(\bm{x})$ be a vector of Koopman eigenfunctions and  $\mathcal{J}(\bm{\varphi}(\bm{x}))$ be its Jacobian matrix. In addition, let $V$ be defined as in \eqref{eq:matModeDecomp}. Then, $P(\bm{x})$ is a right eigenvector of the matrix $V\cdot \mathcal{J}(\bm{\varphi}(\bm{x}))$ with eigenvalue one.
\end{proposition}
\begin{proof}
The time derivative of Eq. \eqref{eq:matModeDecomp} is given by,
\begin{equation}
\begin{split}
    \frac{d}{dt}\bm{x}&=V\frac{d}{dt}\bm{\varphi}(\bm{x})\\
    P(\bm{x})&=V\mathcal{J}(\bm{\varphi}(\bm{x}))P(\bm{x}).
\end{split}
\end{equation}
\end{proof}

\begin{idoexample}[Nonlinear system]
Given the following system,
\begin{equation}\label{eq:exmNonlinearMode}
    \frac{d}{dt}\begin{bmatrix}
    x_1\\x_2
    \end{bmatrix} = \begin{bmatrix}
    x_1\\x_2-x_1^2\end{bmatrix},\quad \begin{bmatrix}
    x_1(0)\\x_2(0)
    \end{bmatrix}=\begin{bmatrix}
    1\\1
    \end{bmatrix}.
\end{equation}
The solution is,
\begin{equation}
\begin{split}
\begin{bmatrix}
    x_1\\x_2
    \end{bmatrix} = \begin{bmatrix}
    1&0\\
    2&-1\end{bmatrix} \begin{bmatrix}
    e^t\\e^{2t}
    \end{bmatrix}.
\end{split}
\end{equation}
The time state-space mappings are,
\begin{equation}
    \begin{split}
        t&=\ln(x_1),\\
        t&=\frac{1}{2}\ln(2x_1-x_2).
    \end{split}
\end{equation}
By choosing $\alpha=1$, $\beta=0$ ,the Koopman eigenfunctions, following \eqref{eq:proofLemTM}, are,
\begin{equation}
    \begin{split}
        \varphi_1(\bm{x})&=x_1,\\
        \varphi_2(\bm{x})&=\sqrt{2x_1-x_2}.
    \end{split}
\end{equation}
The state-space, $\begin{bmatrix}
x_1&x_2
\end{bmatrix}^T$, can be reconstructed by these eigenfunctions as, 
\begin{equation}
    \begin{bmatrix}
    x_1\\x_2
    \end{bmatrix} = \begin{bmatrix}
    \varphi_1(\bm{x})\\2\varphi_1(\bm{x})-\varphi_2(\bm{x})^2\end{bmatrix}=\begin{bmatrix}
    1&0\\
    2&-1\end{bmatrix}\begin{bmatrix}\varphi_1(\bm{x})\\
    \varphi_2(\bm{x})^2
    \end{bmatrix}=V\bm{\varphi}(\bm{x}).
\end{equation}
We observe there are two modes, $[1,2]^T$ and $[0,-1]^T$, which evolve linearly under the nonlinear system \eqref{eq:exmNonlinearMode}. In addition, $V\mathcal{J}(\bm{\varphi}(\bm{x}))=I_{2\times 2}$, hence, $P(\bm{x})$ is an eigenvector for any $\bm{x}$.
\end{idoexample}

We have shown above the strong relation between time state-space mapping and Koopman eigenfunctions. The following proposition states a limitation between the two notions.

\begin{proposition}[Existence of Koopman eigenfunctions with no time state-space mapping]
The state-space mapping is not a necessary condition for the existence of Koopman eigenfunctions. 
\end{proposition}
\begin{proof}
This can be shown by the following simple example. Let us consider the linear system,
\begin{equation}
    \frac{d}{dt}\bm{x}=A\bm{x},\quad \bm{x}(t=0)=\bm{x}_0,
\end{equation}
where $A$ is an $N \times N$ matrix. For simplicity, we assume the eigenvalues, $\{\lambda_i\}_{i=1}^N$, are unique and the eigenvector set, $\{\bm{v}_i\}_{i=1}^N$, is orthonormal. Then, the solution of this system of equations  can be written as,
\begin{equation}
    \bm{x}(t)=\sum_{i=1}^Nb_i \bm{v}_i e^{\lambda_i t},
\end{equation}
where the vector $\bm{b}=\begin{bmatrix}b_1&\cdots&b_N\end{bmatrix}^T$ is chosen according to the initial condition. To form the Koopman eigenfunctions and, correspondingly, the Koopman mode, one should formulate the time state-space mapping. For each eigenvector and eigenvalue of $A$ there is a mapping, expressed as,
\begin{equation}\label{eq:EFlinearDyna}
    t_i(x)=\frac{1}{\lambda_i}\ln\left(\frac{\bm{v}_i^T\bm{x}}{b_i}\right).
\end{equation}
Thus, the Koopman eigenfunctions are,
\begin{equation}\label{eq:KEFlinearDyna}
    \varphi_i(\bm{x})= e^{t_i(\bm{x})}=\left( \frac{\bm{v}_i^T\bm{x}}{b_i}\right)^{\frac{1}{\lambda_i}}.
\end{equation}
This expression can be simplified by applying Def. \ref{def:koopmanFamily},  yielding the following system,
\begin{equation}\label{eq:linearDynaKoop}
    \frac{d}{dt}\begin{bmatrix}\bm{v}_1^T\bm{x}\\
    \vdots\\
    \bm{v}_n^T\bm{x}\end{bmatrix}=\begin{bmatrix}\lambda_1&&\\
    &\ddots&\\
    &&\lambda_n\end{bmatrix}\begin{bmatrix}\bm{v}_1^T\bm{x}\\
    \vdots\\
    \bm{v}_N^T\bm{x}\end{bmatrix}.
\end{equation}
Note that if the eigenvector, $\bm{v}_i$, is complex then the time state-space mapping, Eq. \eqref{eq:EFlinearDyna}, does not exist since it is not well defined. In this case, to create a time state-space mapping, we have to choose one branch from the $\ln$ function. However, the Koopman eigenfunction, Eq. \eqref{eq:KEFlinearDyna}, has a unique value since the exponent cancels the ambiguity of the $\ln$ function. It shows that a Koopman eigenfunction can exist in cases where the  time state-space mapping does not. 
\end{proof}

\section{Koopman Theory for PDE}\label{sec:conKoopman} 
Let us generalize the results above to the continuous setting of Koopman theory, following \cite{nakao2020spectral}. We consider the solution of Eq. \eqref{eq:PDE}, based on the following assumptions.

\begin{assumption}[Proper Operator]\label{assu:properDynamics}
The operator $\mathcal{P}(f(x))$ in Eq. \eqref{eq:PDE} is proper.
\end{assumption}

\begin{lemma}[Continuous $u$]\label{lemma:ConPDESolu}
If the operator $\mathcal{P}$ in Eq. \eqref{eq:PDE} admits Assumption \ref{assu:properDynamics} then the solution is continuous in $t$. 
\end{lemma}
This is quite standard in the theory of PDEs. Basically, letting $u(x,t)$ to be the solution of Eq. \eqref{eq:PDE}, we can write a first order Taylor expansion for the variable $t$  as,
\begin{equation}
    u(x,t+dt)=u(x,t)+\mathcal{P}(u(x,t))\cdot dt + o(dt).
\end{equation}
Since the value of $\mathcal{P}(u(x,t))$ is finite, we get $\abs{u(x,t+dt)-u(x,t)}\to 0$ as $dt\to 0$.

\begin{assumption}[Fr\'echet Differentiability]\label{assu:FrechDiff}
The operator $\mathcal{P}$ is Fr\'echet differentiable a.e. in $\mathcal{H}$.
\end{assumption}

If $\mathcal{P}$ admits Assumption \ref{assu:FrechDiff} then the solution $u(x,t)$ is in $C^1$ a.e. with respect to $t$ (see e.g. \cite{venturi2021spectral}).

\begin{definition}[Time mapping]\label{def:inverSolutionPDE}
    Let $u(x,t)$ be the solution of the dynamical system \eqref{eq:PDE} where $t\in I$. Let  $\Xi(u)$ ($\Xi$ is capital $\xi$) be a functional mapping from the solution $u$ to $t$, i.e.
    \begin{equation}\label{eq:inverSolutionPDE}
        t= \Xi(u).
    \end{equation}
\end{definition}

\begin{lemma}[Differentiation of time mapping ]\label{lemma:timemapDiff} Let the Assumptions \ref{assu:properDynamics} and \ref{assu:FrechDiff} hold. If the time mapping, $t=\Xi(u)$, exists then it admits the following,
\begin{equation}\label{eq:lemmaTMD}
    \inp{\partial \Xi(u(x,t))}{\mathcal{P}(u(x,t))}=1\quad a.e.\,\, \textrm{in}\,\, t\in I.
\end{equation}
\end{lemma}
\begin{proof}
The mapping  $\Xi(u(x))$ is in $C^1$ a.e. in $t\in I$ since $u(x,t)\in C^1,\, a.e.$ with respect to $t$ in $I$. Based on the Brezis chain rule, the time derivative of the mapping is,
\begin{equation}
    1=\frac{d}{dt}t=\frac{d}{dt}\Xi(u)=\inp{\partial \Xi(u(x,t))}{\frac{d}{dt}u(x,t)}=\inp{\partial \Xi(u(x,t))}{\mathcal{P}(u(x,t))}.
\end{equation}
And this expression is valid almost everywhere.
\end{proof}

\begin{proposition}[Condition for the inexistence of a Koopman eigenfunctional]\label{prop:equilibriumPointPDE}
If there is an equilibrium point in $I$ then a nontrivial Koopman eigenfunctional does not exist.
\end{proposition}
\begin{proof}
Let $t_0\in I$ be an equilibrium point and $\phi(u(x,t))$ be a Koopman eigenfunctional. Then, $u(x,t)=const,\,\forall t\in[t_0,b]$. Therefore, Eq. \eqref{eq:eigenFunDer} does not hold for nontrivial $\phi$ for any $\lambda \ne 0$.
\end{proof}

\paragraph{Remark on dynamics with finite time support}
Remark \ref{rem:inexist} is valid also for dynamics of the form of Eq. \eqref{eq:PDE}. Namely, if there exits a time point, $T_{ext}\in I$, for which $\mathcal{P}(u(x,t))=0, \forall t>T_{ext}$, then there is no Koopman eigenfunctional for this dynamics.

\begin{lemma}[Koopman eigenfunctionals induced by a time state-space mapping]\label{lem:mapping2KEnal}
Let the Assumptions \ref{assu:properDynamics} and \ref{assu:FrechDiff} hold and $u(x,t)$ be the solution of Eq. \eqref{eq:PDE}. If there exists a time mapping, $t=\Xi(u)$, then a Koopman eigenfunctional exists.
\end{lemma}
\begin{proof} Given the mapping, $t=\Xi(u)$, we define the following functional,
\begin{equation}\label{eq:proofLemTMnal}
    \phi(u)=e^{\alpha \Xi(u)+\beta}.
\end{equation}
The time derivative of this functional is,
\begin{equation}\label{eq:timeStateDetFunctional}
\begin{split}
    \frac{d}{dt}\phi(u(x,t))&=\frac{d}{d\Xi} e^{\alpha\cdot \Xi(u(x,t))+\beta}\frac{d}{dt}\Xi(u(x,t))=\alpha \phi(u(x,t)) \inp{\partial \Xi(u(x,t))}{\frac{d}{dt}u(x,t)}\\
    &=\alpha \phi(u(x,t)) \inp{\partial \Xi(u(x,t))}{\mathcal{P}(u(x,t))}.
\end{split}
\end{equation}
According to Lemma \ref{lemma:timemapDiff}, $\inp{\partial \Xi(u(x,t))}{\mathcal{P}(u(x,t))}=1$ a.e.. Thus, the function in Eq. \eqref{eq:proofLemTMnal} admits Eq. \eqref{eq:KEFal} for any value of $\beta$, where the corresponding eigenvalue is $\lambda = \alpha$. 
\end{proof}

\begin{theorem}[Sufficient condition for the existence of a Koopman eigenfunctional] \label{theo:sufficientConditionKEFnal} Let the Assumptions \ref{assu:properDynamics} and \ref{assu:FrechDiff}  hold, let $u(x,t)$ be the solution of Eq. \eqref{eq:PDE}, and let there be a real function $f:I\to L$, for which $u(f(t),t)$ is monotonic with respect to $t$. Then, Koopman eigenfunctionals exist in the time interval $I$.
\end{theorem}
\begin{proof}
Let us define the following monotonic function,
\begin{equation}
    g(t)=\int_0^L u(x,t)\delta(x-f(t))dx,
\end{equation} 
where $\delta$ is the Dirac delta. Then, the time mapping is
\begin{equation}
    t=\Xi(u)=g^{-1}\left(\int_0^L u(x,t)\delta(x-f(t))dx\right).
\end{equation}
According to Lemma \ref{lem:mapping2KEnal} there exits a eigenfunctional, which can be expressed by Eq. \eqref{eq:proofLemTMnal}.
\end{proof}

\section{Mode Decomposition based on Time State-Space Mapping} \label{sec:ourDMD} 
\subsection{Bridging between nonlinear spectral decomposition and \ac{KMD}}
Let us recall the dynamical system and its suggested form of solution. We consider the following PDE,
\begin{equation}
    u_t=P(u),
\end{equation}
where $P$ is a nonlinear operator, $u(t=0)=f$.  The solution of this PDE is approximated as
\begin{equation}\label{eq:spatioTempo}
    u(x,t) \approx \sum_{i=1}^m X_i(x)T_i(t).
\end{equation}

We would like to mention two principal PDEs for which this approximation is precise (reaches equality). The first one is linear diffusion and the second is \ac{TV}-flow (see the studies on spectral TV of \cite{gilboa2014total}, \cite{burger2016spectral}, \cite{bungert2019nonlinear}). In both cases, the temporal term $T_i(t)$ are the typical decay profiles of the operator which is dictated by its homogeneity. Whereas the decay profile of linear diffusion is exponential, that of TV-flow is linear. This was generalized by \cite{cohen:hal-01870019,cohen2020introducing}, where it is shown there is a smooth transition between exponential and linear decay for $\gamma$-homogeneous operators, $\gamma \in [0,1))$, see Fig. \ref{Fig:decayProfile}. 
These profiles can be calculated by analyzing an evolution initiated with a single (nonlinear) eigenfunction $f$, admitting $P(f)=\lambda f$. In this case it is simple to check that the evolution is structure preserving. That is, the spatial structure of $f$ is maintained and only its contrast changes throughout the evolution. We thus get a separation of variables and can deduce the time profile. 
It was shown in \cite{bungert2019asymptotic} that the typical decay profile is also the asymptotic behavior of the dynamic (at a time point just before extinction). 

In \cite{gilboa2014total}, \cite{burger2016spectral} it was suggested to perform a decomposition of the signal $f$ by identifying phase transitions of the piecewise linear dynamics of TV, or of gradient flows of one-homogeneous functionals in general. This was performed simply by taking the second time derivative of the flow, where the time-weighted expression $\phi(x,t)=t u_{tt}(x,t)$ was referred to as a spectral component, admitting a simple reconstruction formula, $f=\int_0^\infty \phi(t)dt$. In \cite{gilboa2014total} it was shown that not only the initial condition but the entire solution $u(x,t)$ can be expressed as a weighted integration of the spectral components, 
$$u(x,t)=\int_0^\infty H(t,\tau)\phi(x,\tau) d \tau,$$
where $H(t,\tau)=((\tau-t)/\tau)^+$. Comparing $\phi(x,\tau)$ to $X_i$ and $H(t,\tau)$ to $T_i(t)$ we get an expression similar to \eqref{eq:spatioTempo}, in an integral form. In \cite{burger2016spectral} it was shown that for the discrete one dimensional TV-flow the number of components is finite and we can express the solution $u$ by a sum of weighted spectral components. One can expand the linear decay profile to an infinite some of Koopman eigenfunctions, as done in Eq. \eqref{eq:exmSolutionAsSumKEF}. Hence we can observe that the nonlinear spectral components $\phi$ are actually Koopman modes! These relations and connections are planned to be further investigated in a future work.

When the evolution is TV-flow, the set $\{\phi\}$ is referred to as spectral TV decomposition. In \cite{cohen2020introducing} the idea was generalized to nonlinear decompositions of $\gamma$-homogeneous functionals, $\gamma\in[1,2)$ .
The typical decay profile is a truncated polynomial with fractional degree almost for every value of $\gamma$. Thus, the decomposition was based on fractional calculus, which made this process less accessible numerically.

To bypass the use of fractional calculus it was suggested to apply \ac{DMD} on the gradient descent of the respective homogeneous functional. As discussed earlier, it was shown that recovering the dynamic with \ac{DMD} yields an inherent error, \cite{cohen2021modes}.
A time rescaling method was proposed to improve the \ac{DMD} decomposition. It was shown theoretically that an evolution of a single eigenfunction is constructed accurately and for general signals improvement in the decomposition was achieved. However, a major problem of phase changes in the flow, due to extinction of modes, was not addressed. This is most inherent in flows based on zero-homogeneous operators, common in signal and image processing.
Alternative recent  methods were suggested to improve \ac{DMD}, however none of them tackles well phase transitions in the flow. These methods use machine learning principles in the design of advanced \ac{DMD} algorithms, such as \ac{EDMD} \cite{williams2016extending,williams2015data2,williams2015data} and \acf{KDMD} \cite{kawahara2016dynamic}. Several learning-based approaches suggested to build a data-driven dictionary to reconstruct the dynamics sparsely \cite{bollt2021geometric,li2017extended,pan2021sparsity,rudy2017data}. These works focus on learning the spatial structures that approximate Koopman modes. In other words, these algorithms aim at finding measurements that evolve linearly under the dynamical system.

Since \ac{DMD} is primarily investigated in the context of fluid dynamics, oscillatory flows are more common, and less attention was directed to smoothing or decaying flows, which are most common in image and signal processing. We thus aim at extending the Koopman tools to this type of processes.
System reconstruction based on  finding spatial structures has some limitations, most notably for processes with finitely decaying modes, since the reconstruction of \acp{KEF} may be infinite-dimensional. The reconstruction of a \ac{KEF} as a polynomial of the observation, as in Example \ref{example:FiniteTimeSupport}, contains an infinite vector of measurements, which is highly intractable numerically. 

Our approach is based on the assumption that the observed dynamic has a typical monotone decay profile within a given time interval. Thus, instead of focusing on measurements that decay exponentially, the focus of our algorithm is on finding spatial structures that decay according to a predefined family of profiles.
Let us recall the generalized spectra which was introduced by \cite{Katzir2017Thesis} and \cite{gilboa2018nonlinear}. This work focused on a decomposition induced by the typical decay profile of the respective operator. The spatial structures are deduced from a dictionary containing an overcomplete set of decay profiles. More formally, given a nonlinear dynamic,
\begin{equation*}
    \frac{d\bm{x}}{dt}=P(\bm{x}),
\end{equation*}
with a typical decay profile, $a(t)$, we extract the spatial structure from the solution, $\bm{x}(t)$, with the following optimization problem, 
\begin{equation*}
    \min\arg_{\mathcal{V}}\{\norm{X-\mathcal{V}\mathcal{D}}_\mathcal{F}\},\quad s.t. \min\norm{\mathcal{V}}_0
\end{equation*}
where $X$, $\mathcal{V}$ and $\mathcal{D}$ are defined in Eqs. \eqref{eq:KatzirX}, \eqref{eq:KatzirV}, and \eqref{eq:Katzir}, respectively. In the rest of this section we show that if the decay profile is monotone then the spatial structures resulting from the general spectral decomposition are the Koopman modes of \ac{KMD}.

\subsection{Generalized dynamic mode decomposition}
\paragraph{Spatiotemporal mode decomposition based on a monotone decay profile}
Let us assume the dynamics induces a known typical monotone profile for different spatial structures in the data. The profile, denoted as $a_{\lambda_i}(t)$, varies according to the spatial structure, $\bm{v}_i$, and depends on a parameter $\lambda_i$. In addition, we assume that the solution can be approximate as,
\begin{equation}\label{eq:SMD}
    \bm{x}(t)= \sum_{i=1}^N\bm{v}_i\cdot a_{\lambda_i}(t)+ e
\end{equation}
where $e$ is a small error term. 

Given the time sampling point set $\{t_i\}_0^M$ (not to be confused with time state-space mapping), we define the overcomplete dictionary,
\begin{equation}
    \mathcal{D}=
    \begin{bmatrix}
    a_{\lambda_0}\left(t_0\right)&\cdots&a_{\lambda_0}\left(t_M\right)\\
    &\vdots&\\
    a_{\lambda_L}\left(t_0\right)&\cdots&a_{\lambda_L}\left(t_M\right)\\
    \end{bmatrix},
\end{equation}
where $L$ is large enough. An atom of this dictionary is a row. Since the time profile $a_{\lambda_i}(t)$ is monotone there exists an inverse function for each atom, denoted as,
\begin{equation}\label{eq:inverseAtom}
    t=\xi(a_{\lambda_i}(t)).
\end{equation}
In matrix formulation, for a discrete time setting, this can be written as,
\begin{equation}\label{eq:inverseDictionay}
    \bm{t}=\bm{\xi}(\mathcal{D}),
\end{equation}
where $\bm{t} \in \mathbb{R}^{(L+1)\times M}$. It is assumed that there exists a (sparse) mode matrix $V$ which can approximate the samples of the system $X$ using the dictionary by, 
\begin{equation}\label{eq:ReconDictionary}
    X=V\mathcal{D} + e,
\end{equation}
where $X=\begin{bmatrix}\bm{x}_0&\cdots&\bm{x}_M\end{bmatrix}$ 
and $e$ is a small error term.

\paragraph{Dimensionality Reduction} Following the assumption of \ac{DMD}, we would like to obtain a sparse representation of modes. This problem has been thoroughly investigated and can be formulated as \cite{mairal2014sparse},
\begin{equation}\label{eq:newDMDErrNPJard}
    \begin{split}
        \min_{V}\norm{X-V\mathcal{D}}_\mathcal{F}^2,\quad s.t. \norm{V}_0\leq r,
    \end{split}
\end{equation}
where $\norm{V}_0<r$ indicates the requirement that only up to $r$ columns in $V$ are not zero. This problem is NP-hard and the sparsity constraint is relaxed to solving the following minimization problem,
\begin{equation}\label{eq:newDMDErr}
    \begin{split}
        \min_{V}\norm{X-V\mathcal{D}}_\mathcal{F}^2 +\lambda \norm{V}_1. 
    \end{split}
\end{equation}
The solution of \eqref{eq:newDMDErr} is the minimizer of the left term when the nonzero entries in each mode are at least $\lambda$ (see algorithm 6 p. 153 in \cite{mairal2014sparse}).

In general, there are several well known algorithms to recover the modes when the dictionary is known (see \cite{elad2010sparse}). We note that our problem is somewhat more difficult than the common signal processing case since the atoms in the dictionary are highly coherent (strongly correlated). Here, we apply the implementation from \cite{mairal2014spams} for the Lasso algorithm (Eq. \eqref{eq:newDMDErr}) with a fine-tuning post-processing stage (\ref{appsec:naive}) to solve this problem. The output of this algorithm is $\hat{V}$ and $\mathcal{\hat{D}}$, where each column in the matrix $\hat{V}$ contains a mode and $\mathcal{\hat{D}}$ has the corresponding atoms, taken from the dictionary $\mathcal{D}$. The entire dynamics can be approximated as,
\begin{equation}\label{eq:koopmanModeApprox}
    X\approx \hat{V}\mathcal{\hat{D}},
\end{equation}
where $\approx$ denotes equality in the sense of Eq. \eqref{eq:newDMDErr}.

\paragraph{Approximation of Koopman eigenfunctions} Given the modes $\hat{V}$ and the data matrix $X$ and assuming $\hat{V}^T\hat{V}$ is invertible, one can express the dictionary as
\begin{equation}\label{eq:decayProfileMD}
    \mathcal{\hat{D}}\approx (\hat{V}^T\hat{V})^{-1}\hat{V}^TX.
\end{equation}
This reconstruction of the dictionary is necessary to be in the argument of the time state-space mapping, Eqs. \eqref{eq:inverseAtom} and \eqref{eq:inverseDictionay}, as follows,
\begin{equation}\label{eq:tssmMD}
    \bm{t}=\bm{\xi}(\mathcal{\hat{D}})=\bm{\xi}((\hat{V}^T\hat{V})^{-1}\hat{V}^TX).
\end{equation}
Thus, we can express with the dynamic measurements an exponential function. According to Eq. \eqref{eq:proofLemTM}, the \acp{KEF} are given by,
\begin{equation}\label{eq:KEFMD}
    \bm{\varphi}(X)=e^{\bm{t}(X)}=e^{\bm{\xi}((\hat{V}^T\hat{V})^{-1}\hat{V}^TX)}.
\end{equation}
We summarize this algorithm in Algo. \ref{algo:DMD}.

\begin{algorithm}[phtb!] \caption{Koopman Mode Approximation}
\begin{algorithmic}[1]
		\Inputs{Data sequence $\{{\bm{x}_k}\}_0^{N}$ and typical profile $a_{\lambda}(t)$}.
		\State Find modes $\hat{V}$ and dictionary $\mathcal{\hat{D}}$ (for example invoke Algo \ref{algo:SR}).
		\State Formulate the decay profiles with the modes $\hat{V}$ and the data $X$, Eq. \eqref{eq:decayProfileMD}.
		\State Formulate the time state-space mapping, Eq. \eqref{eq:tssmMD}.
		\Outputs{Extract \acp{KEF} from the observations by Eq. \eqref{eq:KEFMD}.}
    \end{algorithmic}
    \label{algo:DMD}
\end{algorithm}

\paragraph{Relation between spatiotemporal mode decomposition and \ac{KMD}} The definition of \ac{KMD} is to express the state-space vector as spatiotemporal mode decomposition where the temporal terms are exponential functions (\acp{KEF}). This can be done easily by extracting the time variable $t$ from Eq. \eqref{eq:KEFMD} and plugging it in Eq. \eqref{eq:SMD}. Then, the typical decay profile $a_{\lambda_i}(t)$ can be expressed using a Taylor series (under sufficient smoothness conditions). By variation of parameter, the \ac{KMD} is obtained (see Example \ref{exm:ZeroHomoDynamics}). 

Note that the above presentation is only intended to show a possible algorithmic path that is implied by our analysis. We limit the scope of our discussion here and leave for future work important issues, such as spectrum and system reconstruction accuracy, dimensionality reduction, robustness to noise, and prediction capacity, for more details on these concepts see \cite{gavish2014optimal,lu2020prediction}.

\section{Examples}\label{sec:results}

In this section, we apply the theory to a few examples. We examine the following: system reconstruction; global controllability; mode decomposition based on a dictionary of monotone profiles; and finding eigenfunctionals in partial differential equations. 

\begin{idoexample}[System Reconstruction and Global Controllability]
This example is based on \cite{mauroy2020koopman} (p. 10). Given the system,
\begin{equation}\label{eq:sec6example1}
    \frac{d}{dt}x(t)=P(x)+u=x-x^3 + u,
\end{equation}
we would like to obtain global controllability via a Koopman eigenfunction according to Remark \ref{remark:GlobalControllability}. Note that there are three equilibrium points $-1,0$ and $1$ with \ac{ROA}s: $\mathcal{RA}(-1)=(-\infty,0),\,\mathcal{RA}(0)=\{0\}$, and $\mathcal{RA}(1)=(0,\infty)$, respectively. The solution of this equation is,
\begin{equation}
\begin{split}
    t(x) &= \ln\left(\frac{x}{\sqrt[]{1-x^2}}\right)+C .
\end{split}
\end{equation}
According to Theorem \ref{theo:sufficientConditionKEF} one of the Koopman eigenfunctions is,
\begin{equation}
    \varphi(x)=e^{t(x)}=\frac{x}{\sqrt[]{1-x^2}}.
\end{equation}
We set the input $u$ to,
\begin{equation}
    u=-\mathcal{J}(\varphi)^{-1}\varphi+w,
\end{equation}
where $w$ is the input after feedback linearization. The Jacobian matrix is simply the derivative of $\varphi$ with respect to $x$,
\begin{equation}
\begin{split}
    \mathcal{J}(\varphi)&=\left(1-x^2\right)^{-\frac{3}{2}},
\end{split}
\end{equation}
yielding,
\begin{equation}
    u=-\mathcal{J}(\varphi)^{-1}\varphi+w = -\left(1-x^2\right)^{\frac{3}{2}}\frac{x}{\sqrt[]{1-x^2}}+w=-x(1-x)^2+w.
\end{equation}
Substituting this input in the dynamical system, Eq. \eqref{eq:sec6example1}, we get the following,
\begin{equation}
    \frac{d}{dt}x(t)=P(x)+u=x-x^3 + u=x-x^3 -x(1-x)^2+w=w.
\end{equation}
This system is linear and controllable.
\end{idoexample}

\begin{idoexample}[\acl{TV} eigenfunctional]
A very common PDE in image processing is the gradient descent flow with respect to the total-variation (TV) functional \cite{bellettini2002total}, which for smooth functions $u$ can be expressed as,
\begin{equation}
    J_{TV}(u(x))=\inp{\abs{\nabla{u(x)}}}{1}.
\end{equation}
The gradient descent flow for this non-smooth convex functional is defined by,
\begin{equation}
    u_t=\mathcal{P}\in -\partial J_{TV}(u),  \qquad u(t=0)=u_0,
\end{equation}
where $\partial J_{TV}(u)$ denotes the subdifferential of TV at $u$. The flow is known also as the 1-Laplacian flow.
When $x\in\mathbb{R}$ the solution is piece-wise linear, at any time interval $\mathcal{I}_j$ the solution admits, \cite{cohen2021Total}, 
\begin{equation}
    u(x,t) = h_{1,j}(x)+h_{2,j}(x)\lambda_j t.
\end{equation}
In addition, it was shown by \cite{burger2016spectral,cohen2021Total} that the two modes are orthogonal, $h_{1,j}\perp h_{2,j}$. Thus, at each interval there are two eigenfunctionals, the trivial one and the second one, corresponding to the linearly evolving mode,
\begin{equation}
    \phi(u)=e^{\frac{\inp{h_{2,j}}{u}}{\norm{h_{2,j}}^2\lambda_j}}.
\end{equation}
\end{idoexample}

\begin{idoexample}[Nonlinear PDE \# 2] Let the solution of Eq. \eqref{eq:PDE} be,
\begin{equation}\label{eq:dynamicExm}
    u(x,t)=v_1(x)\cdot a_1(t)+v_2(x)\cdot a_2(t).
\end{equation}
The solution $u(x,t)$ and the spatial structures $v_i(x)$, $i\in \{1,2\}$, are depicted in Fig. \ref{fig:refSolution}. The decay profile is of the form of $a_i(t)=(1+\lambda_i t)^+,$ where $\lambda_1=1/10$ and $\lambda_2=1/30$.
\begin{figure}[htbp!]
    \centering
    \subfloat[\bf{$u(x,t)$}]{
    \includegraphics[width=0.3\textwidth,valign=t]{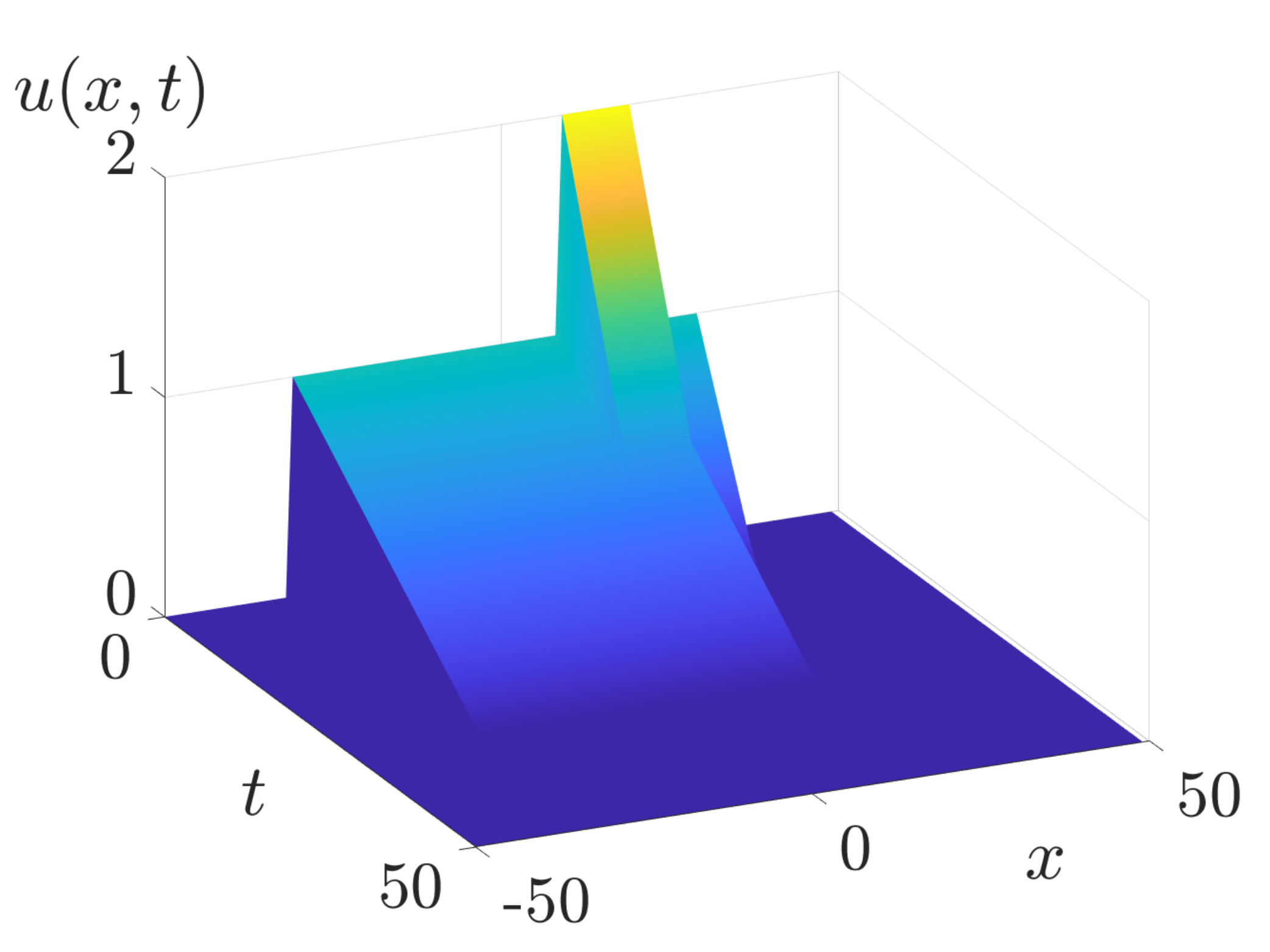}
    \label{subfig:Xall}}
    \subfloat[\bf{$v_1$}]{
    \includegraphics[width=0.3\textwidth,valign=t]{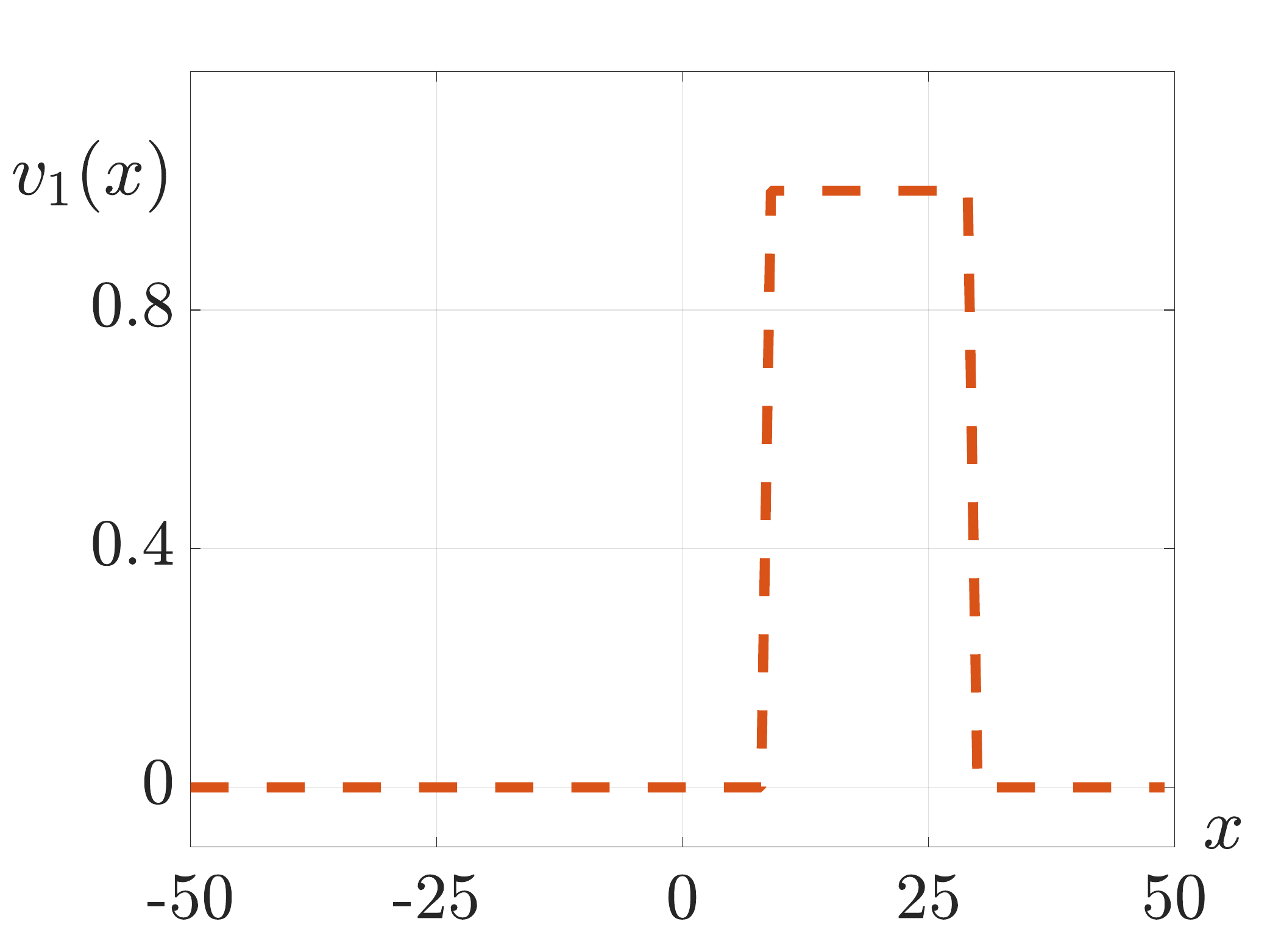}
    \label{subfig:RefMode1}}
    \subfloat[\bf{$v_2$}]{
    \includegraphics[width=0.3\textwidth,valign=t]{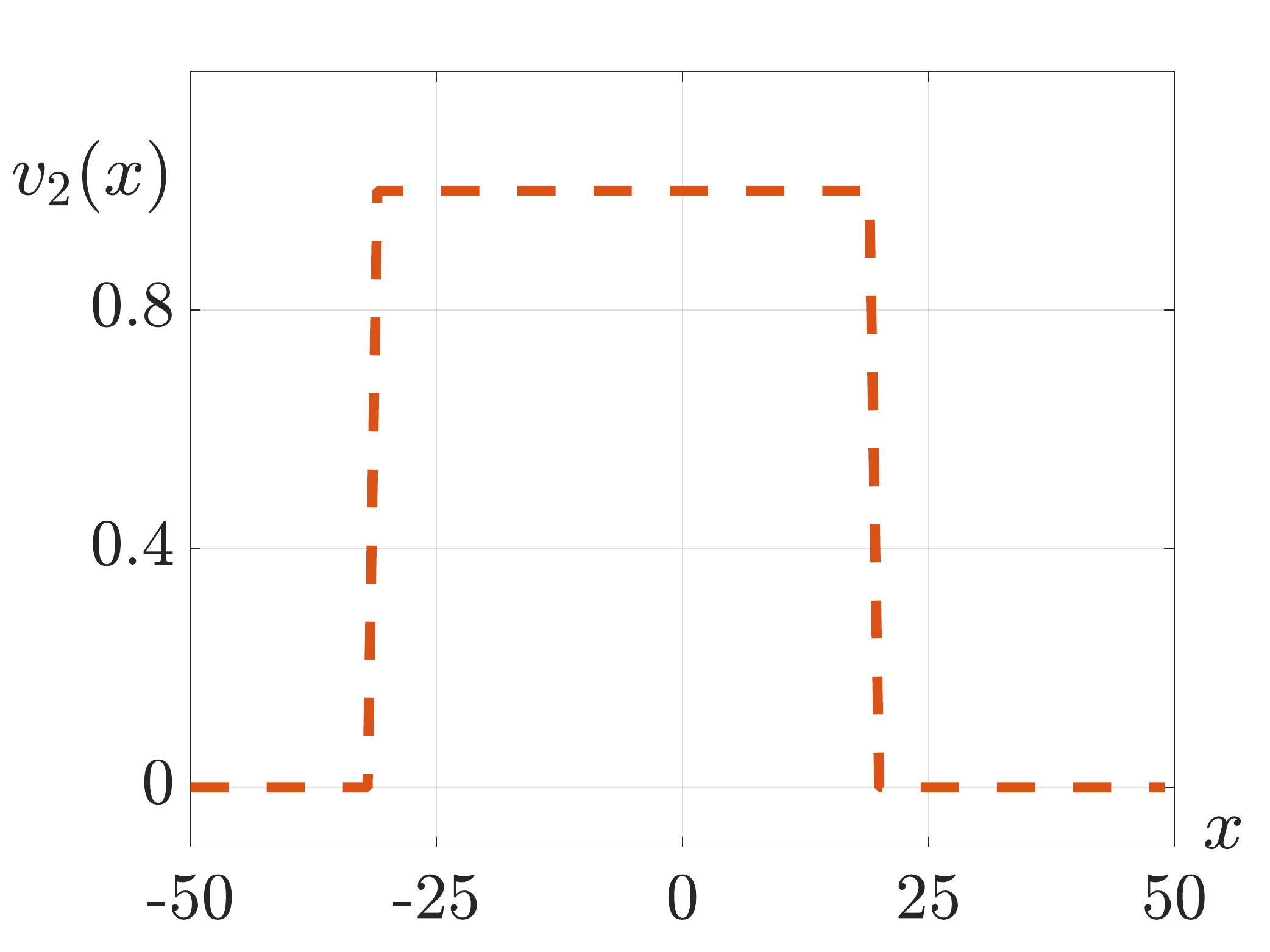}
    \label{subfig:RefMode2}}
\caption{(a) The solution of Eq. \eqref{eq:dynamicExm}. On the right, the spatial structures (modes), $v_1$ (b) and $v_2$ (c). They evolve with linear decay at a rate of $\lambda_1=1/10$ and $\lambda_2=1/30$, respectively.}
\label{fig:refSolution}
\end{figure}
\ac{DMD} yields the decomposition depicted in Fig. \ref{fig:modesDMDLinear}. The modes are complex and each of them is depicted in two graphs, the real and the imaginary parts (Fig. \ref{subfig:modesDMDLinear}). It demonstrates the limitations of \ac{DMD} in systems with typical dynamics which are not exponential.

\begin{figure}[htbp!]
    \centering 
    \captionsetup[subfigure]{justification=centering}
    \subfloat[{\bf{\acl{DMD}.}} Two plots on the left: real and imaginary values of the first \ac{DMD} mode, compared to $v_1$ (dashed). Two plots on the right: real and imaginary values of the second \ac{DMD} mode, compared to $v_2$ (dashed).]{
    \includegraphics[trim=220 0 250 0, clip,width=0.975\textwidth]{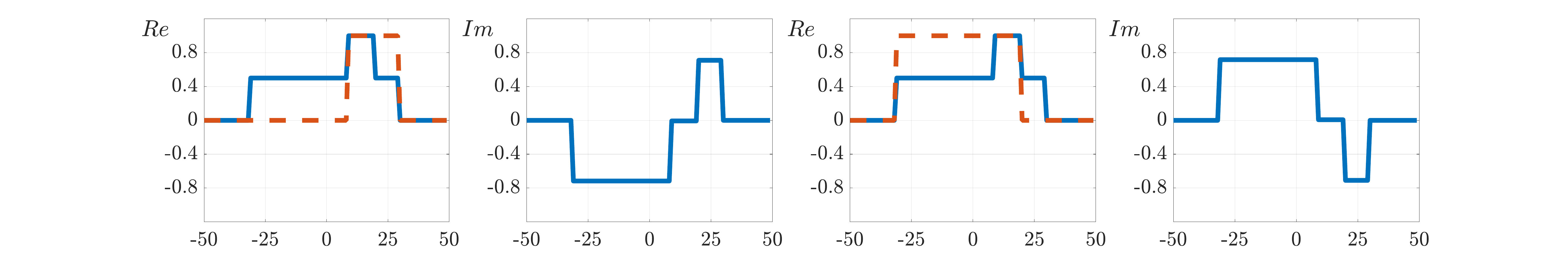}
    \label{subfig:modesDMDLinear}}\\
    \subfloat[{\bf{Reconstruction}}]{
    \includegraphics[width=0.4\textwidth,valign=c]{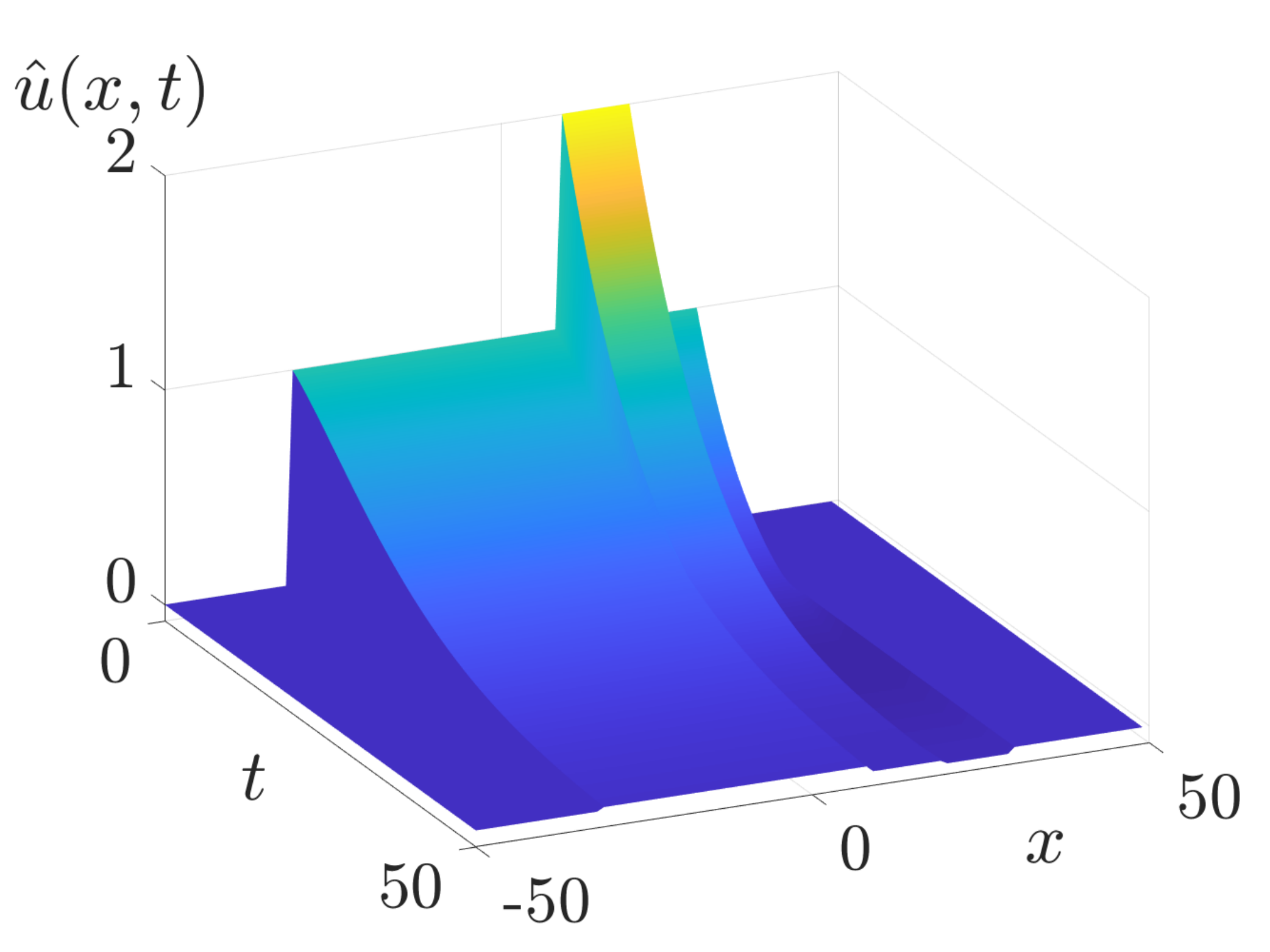}
    \label{subfig:Xkova}}
    \subfloat[\bf{Error $u(x,t)-\hat{u}(x,t)$}]{
    \includegraphics[width=0.4\textwidth,valign=c]{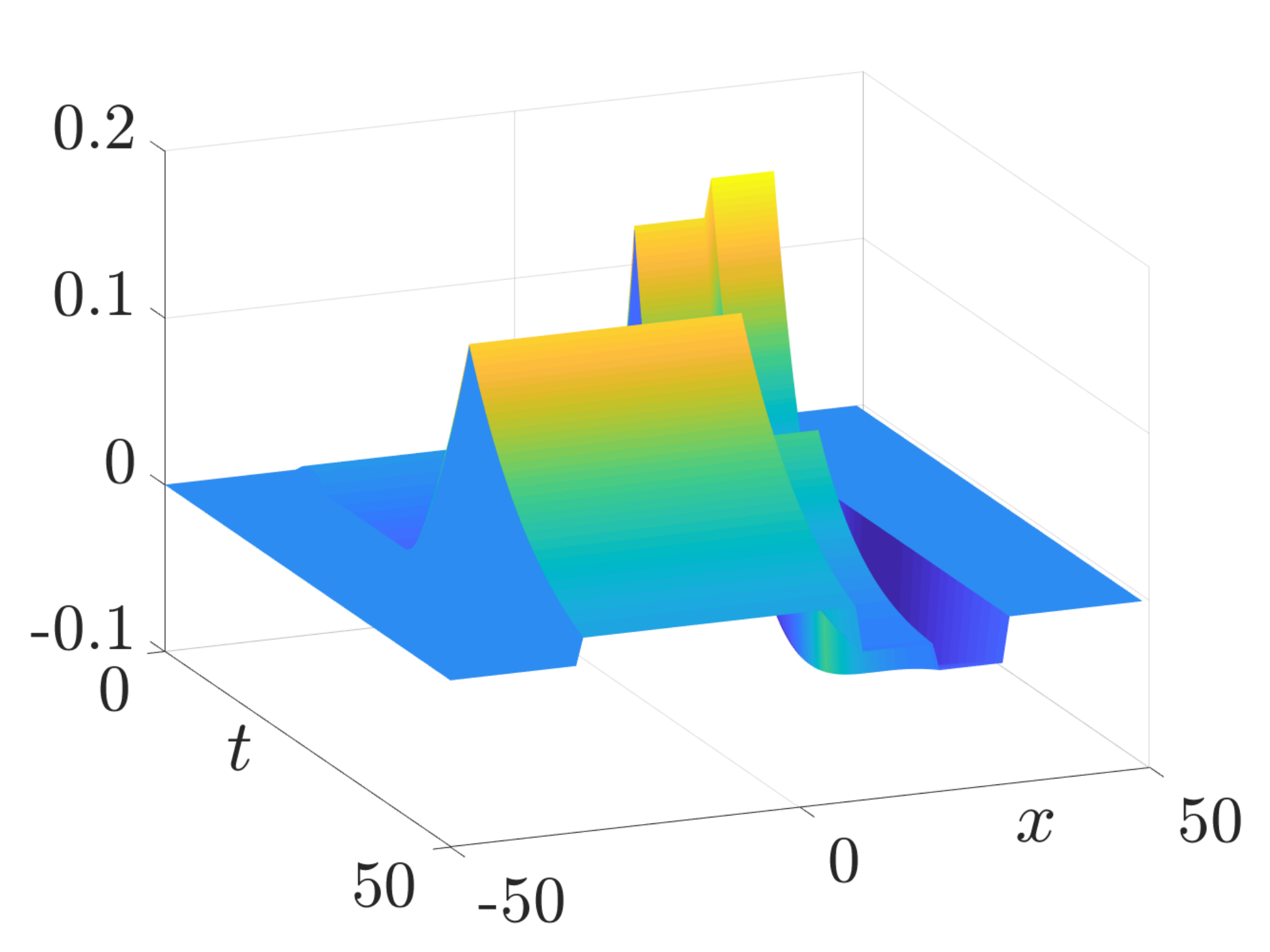}
    \label{subfig:Xdiff}}
    \caption{{\bf{Dynamic Mode Decomposition and Reconstruction.}} Top row: First two \ac{DMD} modes, compared to $v_i$. Bottom row: Reconstruction through Eq. \eqref{eq:dynamicApproDiscrete} (left) and the corresponding error (right). We can observe the dynamics is not reconstructed well and the error is significant.}
    \label{fig:modesDMDLinear}
\end{figure}

The decomposition resulted from Algo. \ref{algo:SR} is depicted in Fig.  \ref{Fig:errorRecon}. The modes are shown in Fig. \ref{subfig:sparseLinear} and recover the modes accurately. The entire dynamics reconstruction is given in Fig \ref{subfig:XkovaSparse} with the corresponding error in Fig. \ref{subfig:XdiffSparse}.
\begin{figure}[htbp!]
    \centering 
    \captionsetup[subfigure]{justification=centering}
    \subfloat[{\bf{Sparse Mode Decomposition}} - The modes resulting from Algo. \ref{algo:SR} (blue) and the actual modes of the dynamics (dashed red).]{
    \includegraphics[trim=0 0 0 0, clip,width=0.667\textwidth]{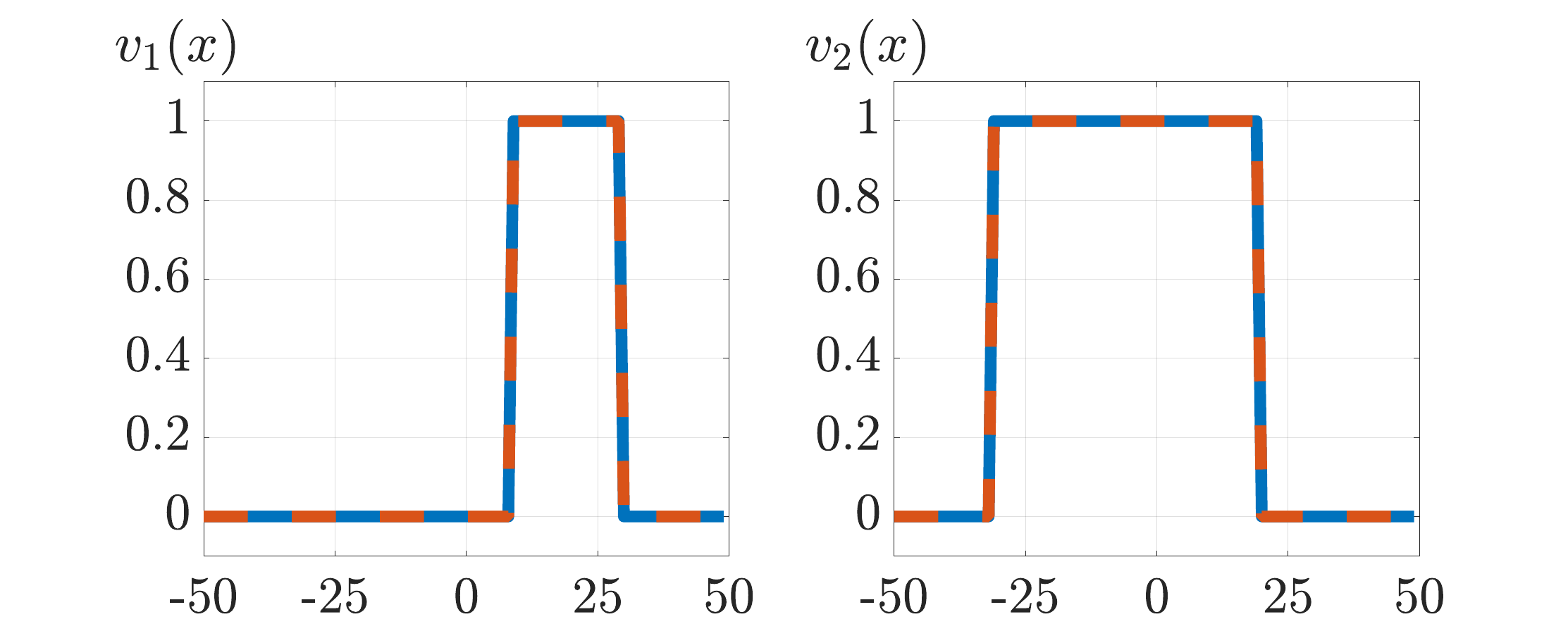}
    \label{subfig:sparseLinear}}\\
\subfloat[\bf{Reconstruction}]
  {
\includegraphics[width=0.4\textwidth,valign=c]{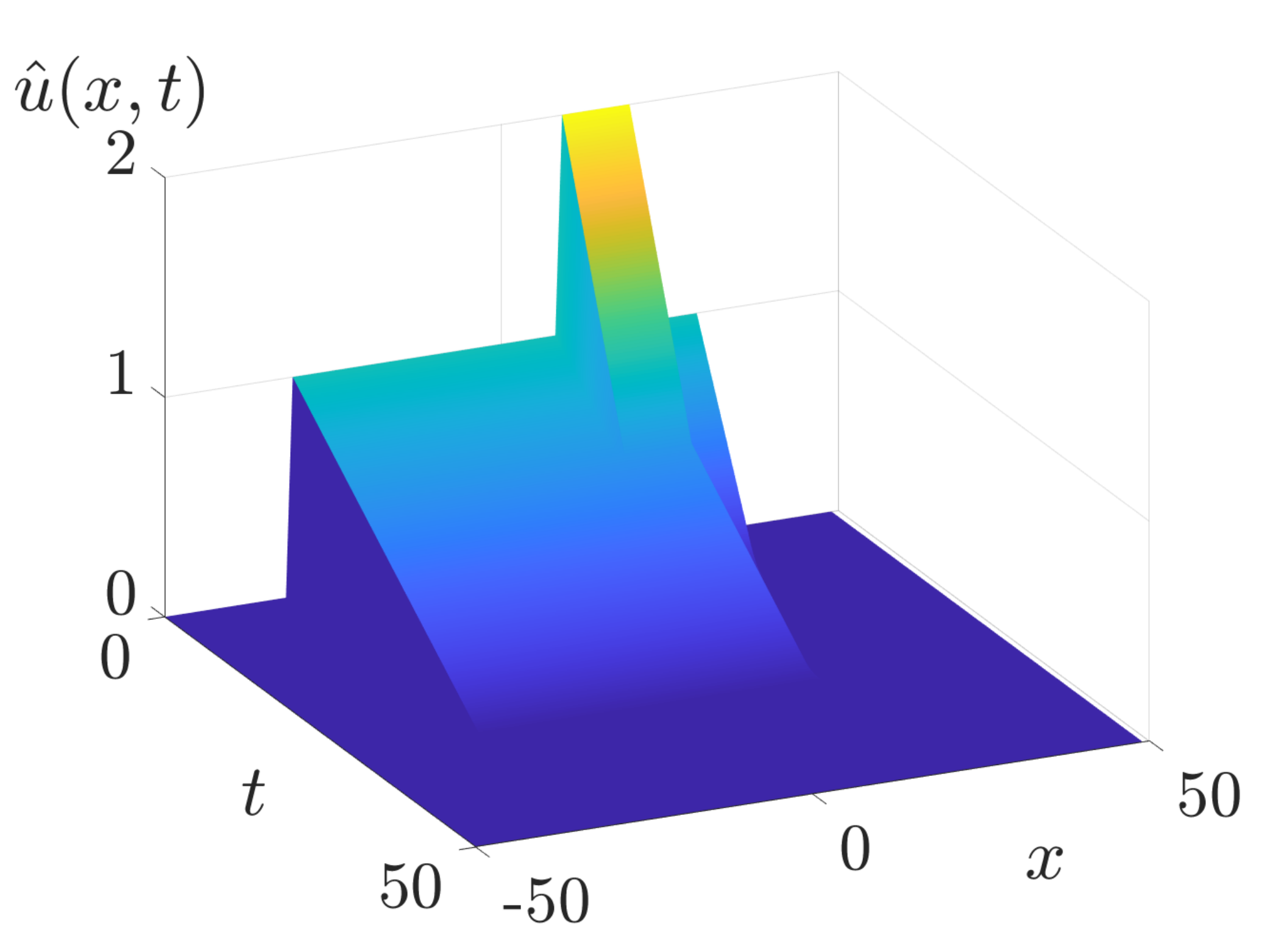}
\label{subfig:XkovaSparse}
}
\subfloat[\bf{Error $u(x,t)-\hat{u}(x,t)$}]
  {
\includegraphics[width=0.4\textwidth,valign=c]{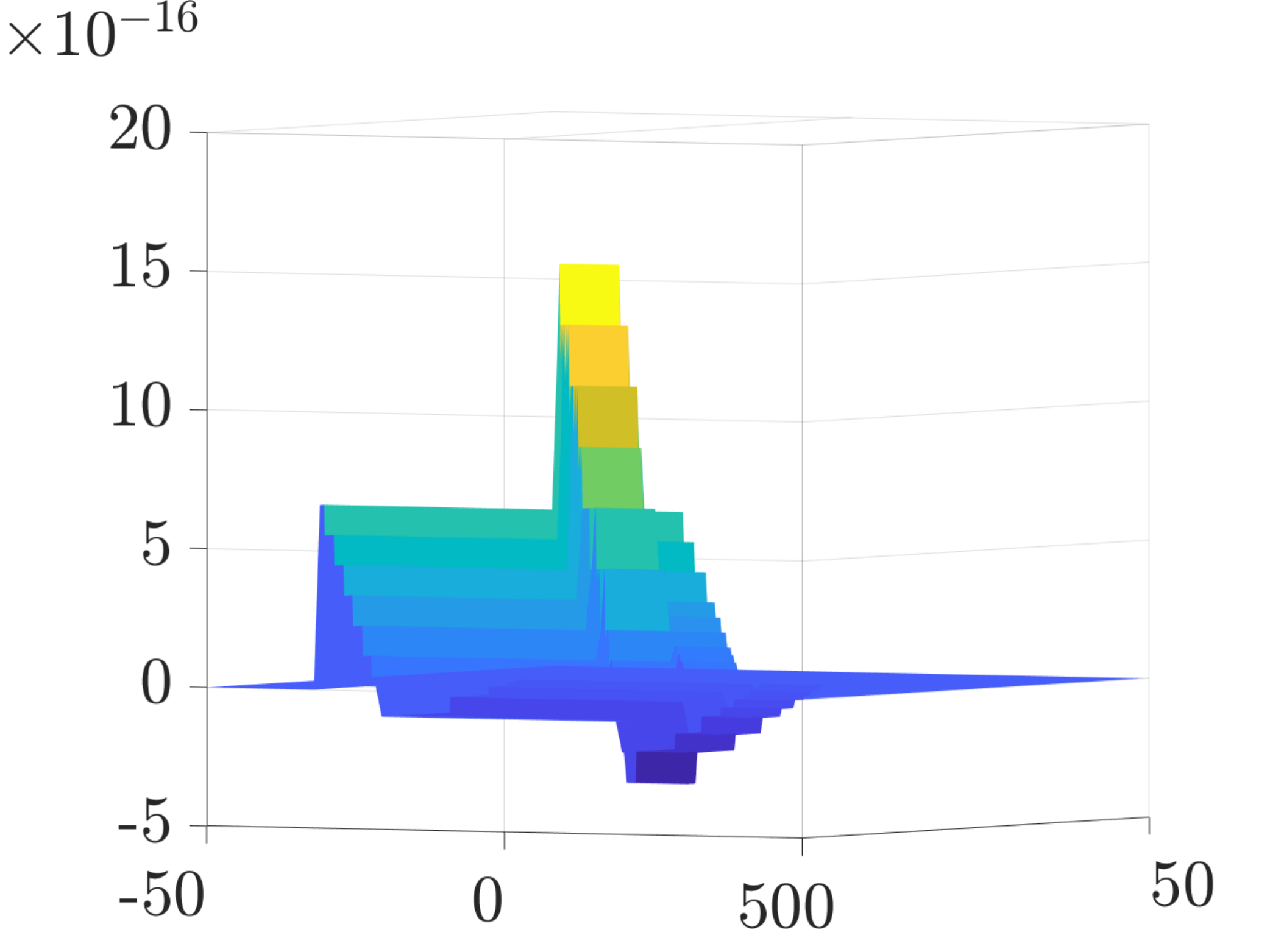}
\label{subfig:XdiffSparse}
}
\caption{{\bf{Sparse Mode Decomposition and Reconstruction (Algo. \ref{algo:SR})}} - (a) Sparse mode decomposition compared to the modes, $v_1$ and $v_2$. (b) Dynamic reconstruction with Algo. \ref{algo:SR} (Eq. \eqref{eq:koopmanModeApprox})
(c) The corresponding error. Correct modes are obtained yielding close to perfect reconstruction of the dynamics.}
\label{Fig:errorRecon}
\end{figure}
Having the modes, we can find the eigenfunctionals,
\begin{equation}\label{eq:examEFal}
    \bm{\phi}(t)=\begin{bmatrix}\inp{v_1}{v_1}&\inp{v_1}{v_2}\\
    \inp{v_2}{v_1}&\inp{v_2}{v_2}\end{bmatrix}^{-1}\begin{bmatrix}\inp{v_1(x)}{u(x,t)}\\\inp{v_2(x)}{u(x,t)}\end{bmatrix}-\begin{bmatrix}1\\1\end{bmatrix}.
\end{equation}
They are depicted in Fig. \ref{Fig:EFnal}. One can see that the eigenfunctionals are valid until the vanishing points. The first mode vanishes at $t=10$ and the second at $t=30$.

\begin{figure}[phtb!]
\centering
\captionsetup[subfigure]{justification=centering}
\subfloat[Eigenfunctional \#1]
  {
\includegraphics[width=0.45\textwidth,valign=c]{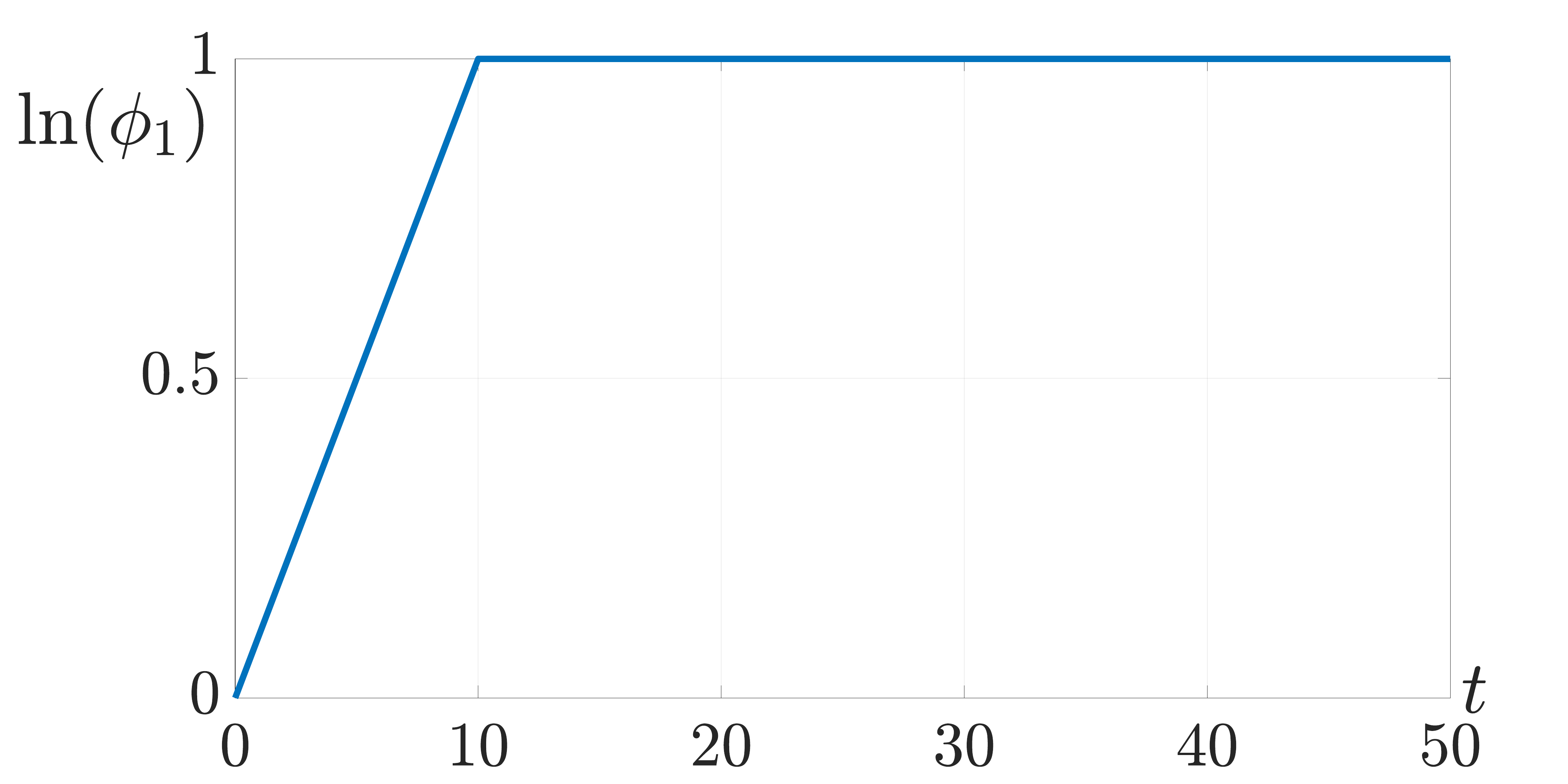}
\label{subfig:J1}
}
\subfloat[Eigenfunctional \#2]
  {
\includegraphics[width=0.45\textwidth,valign=c]{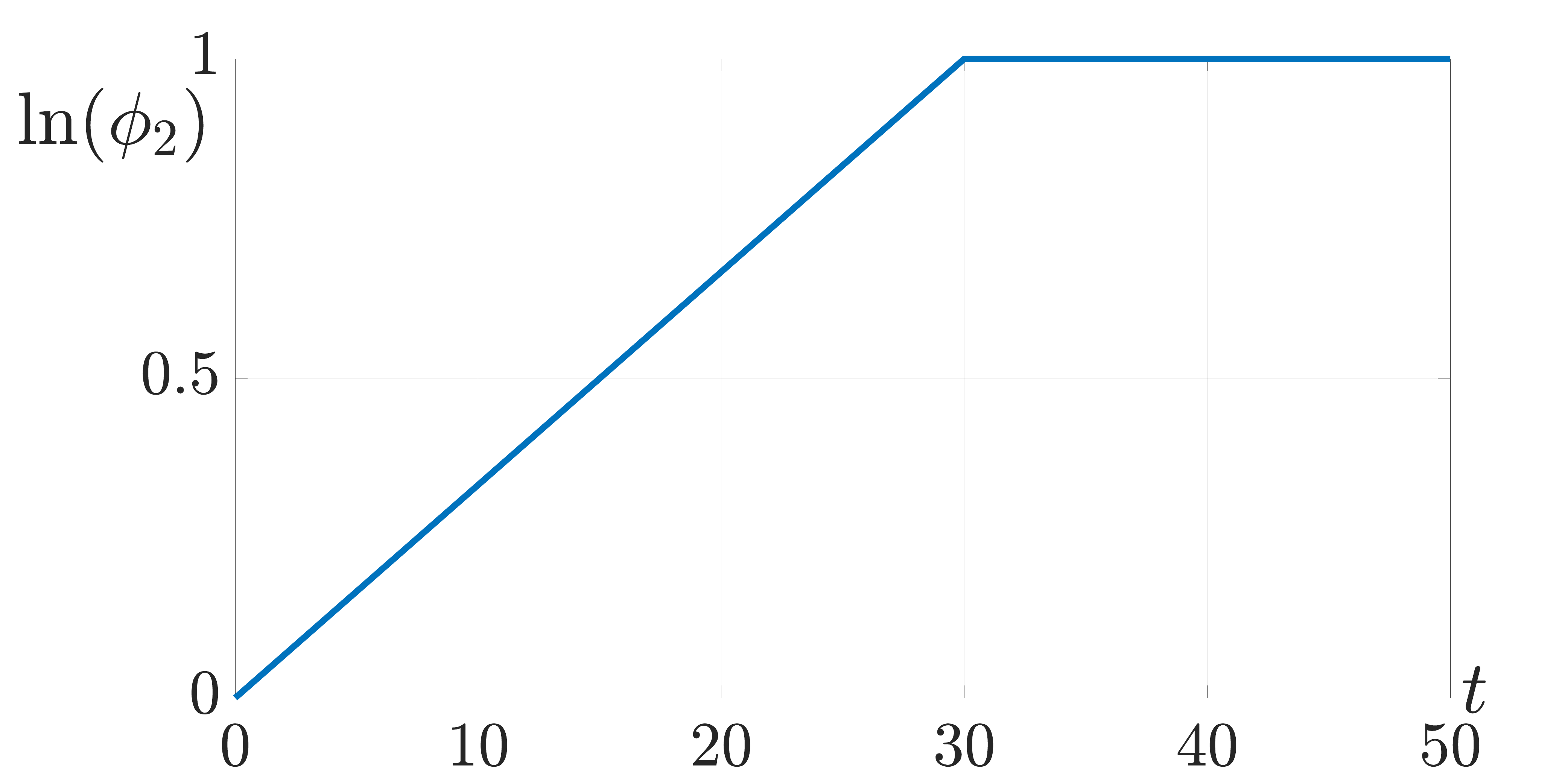}
\label{subfig:J2}
}
\caption{{\bf{Eigenfuntionals}} - based on the monotone decay profile dictionary. These are the eigenfunctionals formulated in Eq. \eqref{eq:examEFal}.   }
\label{Fig:EFnal}
\end{figure}

\end{idoexample}

\section{Conclusion}
This work investigates ways to broaden the use of tools from Koopman theory for the analysis of local and nonlocal PDE's emerging in image and signal processing. We focus on evolution of smoothing processes with possible phase transitions in the dynamics, inherent in zero-homogeneous operators.
We discuss necessary and sufficient conditions for the existence of Koopman eigenfunctions. We examine \ac{KMD}, system reconstruction, global linearity, controllability, and observability through Koopman theory. These insights highlight some limitations of \ac{DMD}. With the technique of time state-space mapping, we show how  conservation laws emerge naturally from any \ac{KEF}. In addition, we justify the approximation of \ac{EDMD} based on this mapping.

The classical \ac{DMD} accurately evaluates \ac{KMD} as long as \acp{KEF} are linear combinations of the observations and \ac{KMD} is finite-dimensional. However, \ac{DMD} has clear limitations in four different settings:
1) The typical decay profile of the system is not exponential; 2) One Koopman mode is associated with multiple eigenvalues; 3) There is an equilibrium point in the time interval $I$; 4) Koopman modes do not exist for all $t$ in $I$. Another limitation  emerges when the dynamic $P$ is in $C^0$ almost everywhere. In this case, some of the modes might vanish at different times, as we see in the total-variation  flow.

We suggest a new type of decomposition to overcome these  fundamental problems. It is based on inverse time state-space mapping of injective curves. We implement this method using overcomplete dictionaries of monotone profiles, typical to the dynamics. This decomposition coincides with a basic assumption of \ac{DMD} -- a flow can be sparsely represented by a few dominant modes. We show our decomposition yields Koopman modes. This work can lead to many interesting connections between decomposition, signal representation, nonlinear PDE's and their relation to Koompan theory.

{\bf{List of Symbols}}
\addcontentsline{toc}{section}{List of Symbols}
\begin{longtable}{lp{0.6\textwidth}}
$\bm{x}_i$ & The $i$ the sample of the state vector belongs to $\mathbb{R}^N$\\
$X$ & Contains the samples of the dynamics $X=\begin{bmatrix}\bm{x}_0&\cdots&\bm{x}_M\end{bmatrix}$ belongs to $\mathbb{R}^{N\times (M+1)}$\\
$\mathcal{U}$&A matrix where $\mathcal{U}_{i,j}=u(x_i,t_j)$\\
$H$&A auto-correlation matrix of the set $\{h_i\}_{i=1}^M$\\
$\bm{Hu}$&A vector where $\bm{Hu}_i=\inp{h_i(x)}{u(x,t)}$\\
$\bm{h}(x)$&A vector $\bm{h}_i(x)=h_i(x)$\\
$\mathcal{V}$& Contains the main spatial structures $\{\bm{v}_i\}$\\
$\mathcal{D}$& A dictionary of a family of a decay profile\\
$P$& A (nonlinear) function $P:\mathbb{R}^N\to \mathbb{R}^N$ in $C^1$ a.e.\\
$t$& Time index where $t\in \mathbb{R}^+$\\
$g$& This is an observation function of the state vector $\bm{x}$,  $g:\mathbb{R}^N \to \mathbb{R}$\\
$K_P^\tau$& The Koopman operator. The superscript denotes the time parameter and the subscript denotes the dynamical system\\
$I$& An interval $[a,b]$ in the time axis\\
$\varphi(\bm{x}(t))$& A Koopman eigenfunction\\
$\lambda$& A Koopman eigenvalue\\
$\nabla$&The gradient of a function\\
$^T$& denotes the transform\\
$\mathcal{H}$& A Hilbert space\\
$\mathcal{P}$& An (nonlinear) operator $\mathcal{P}:\mathcal{H}\to\mathcal{H}$\\
$Q$& A (nonlinear) proper, lower-semicontinuous functional $Q:\mathcal{H}\to\mathbb{R}$\\
$\phi(\cdot)$& A Koopman eigenfunctional\\
$\bm{v}_i$& A preserved spatial shape under the dynamics $P$\\
$h_i(x)$& A preserved spatial shape under the dynamics $\mathcal{P}$\\
$a_i(t)$& The time profile corresponding to the $i$th preserved spatial shape\\
$\gamma$, $\gamma-1$&Denote the homogeneity degrees of a functional and its variational derivative, respectively.\\
$X_0^{M-1},X_1^M$& Data matrices $[\bm{x}_0,\cdots,\bm{x}_{M-1}],[\bm{x}_1,\cdots,\bm{x}_{M}]$\\
$U,\Sigma, V$& \acf{SVD} of $\bm{x}_0^{N-1}$\\
$U_r, V_r$& Sub-matrices of $U, V$ containing the first $r$ columns\\
$\Sigma_r$&Sub-matrix of $\Sigma$ containing the most significant $r$ eigenvalues of the \ac{SVD} which are the diagonal of $\Sigma$\\
$\mathcal{X}$& The curve in $\mathbb{R}^N$ representing the solution $\bm{x}$\\
$\xi(\cdot)$&A mapping from the curve $\bm{x}(t)$ to the time variable $t$\\
$\bm{\varphi}(\bm{x})$& A Koopman mode\\
$\mathcal{J}(\bm{\varphi}(\bm{x}))$&The Jacobian of Koopman mode $\bm{\varphi}(\bm{x})$\\
$\Xi$& A functional mapping from $u(x,t)$ to $t$\\
$\phi$&An eigenfunctional\\

  \end{longtable}

\appendix

\section{The \acl{DMD} steps}\label{appsec:DMDsteps}

\paragraph{Coordinate representation} Given $N$ observations of the dynamical system, Eq. \eqref{eq:disDynamicalSystem}, we form the data matrices as 
\begin{equation}\label{eq:DMDdataForming}
X_0^{M-1}=[\bm{x}_0,\bm{x}_1,\cdots,\bm{x}_{M-1}], \quad X_1^{M}=[\bm{x}_1,\bm{x}_2,\cdots,\bm{x}_{M}]\in \mathbb{R}^{N \times {M}}    
\end{equation}
where $\bm{x}_k=\bm{x}(t_k)$. To find the spatial structures the \ac{SVD} is applied on the data matrix,
\begin{equation}\label{eq:SVD}
    X_0^{M-1} = U\Sigma V^* .
\end{equation}
where $V^*$ is the conjugate transpose of $V$. The columns of $U$ span the column space of $X_0^{M-1}$. Thus, the spatial structures are represented by its coordinates
\begin{equation}
    \bm{c}_k = U^*\bm{x}_k.
\end{equation}

\paragraph{Dimensionality reduction} Assuming the data is embedded in subspace spanned by the first $r$ columns of $U$. Then, the coordinates related to that subspace is
\begin{equation}\label{eq:coordinatesK}
    \bm{c}_{r,k}=U_r^*\cdot \bm{x}_k.
\end{equation}

\paragraph{Linear mapping} Following the second assumption of the \ac{DMD}, there is a linear mapping, $F$, from $\bm{c}_{r,k}$ to $\bm{c}_{r,k+1}$. The linear mapping, $F$, minimizes the \ac{DMD} error, given by
\begin{equation}\label{eq:DMDERR}
    F=\arg\min_F \norm{F\cdot C_{r,0}^{M-1}-C_{r,1}^M}_\mathcal{F}^2,
\end{equation}
where $\norm{\cdot}_\mathcal{F}$ denotes the Frobenius norm and 
\begin{equation}\label{eq:DMDdimeReduction}
\begin{split}
    C_{r,0}^{M-1}&=U_r^*X_0^{M-1},\quad
    C_{r,1}^{M}=U_r^*X_1^M.
\end{split}
\end{equation} 

The linear mapping, $F$, is the optimal linear mapping in the sense of the \ac{DMD} error, Eq. \eqref{eq:coordinatesK}, and we write the coordinate dynamic as
\begin{equation}\label{eq:approxLinearMapping}
    \bm{c}_{r,k+1}\approx F\cdot \bm{c}_{r,k}.
\end{equation}
Then, we can write the dynamic for all $k$ as
\begin{equation}
    \bm{c}_{r,k}\approx F^k\cdot \bm{c}_{r,0}.
\end{equation}

\paragraph{Modes, eigenvalues, and coefficients} Now, we would like to summarize the discussion above and to depict the dynamics as a linear system.
In general, we can reconstruct a sample at step $k$ from its coordinates as
\begin{equation}
\bm{\tilde{x}}_{k}=U_r\bm{c}_{r,k}.    
\end{equation}
In addition, if $F$ is diagonalizable it can be formulated as
\begin{equation}\label{eq:diagonalizableF}
    F=WDW^{-1},
\end{equation}
where $W$ contains the right eigenvectors of $F$, and $D$ is a diagonal matrix whose entries are the eigenvalues of $F$.

Then, the dynamic can be simplified as
\begin{equation}
    \bm{\tilde{x}}_{k}\approx U_r\cdot F^k\cdot U_r^*  \bm{x}_0=U_r\cdot WD^kW^{-1}\cdot U_r^*  \bm{x}_0
\end{equation}

Now, let us define the modes, $\{\phi_i\}_{i=1}^r$, eigenvalues,$\{\mu_i\}_{i=1}^r$, and coefficients, $\{\alpha_i\}_{i=1}^r$. 

\emph{Modes} are defined as $\Phi=\begin{bmatrix}{\bm{\phi}}_1&\cdots&\bm{\phi}_r\end{bmatrix}=U_rW$. 

\emph{Eigenvalues} are the diagonal entries of the matrix $D$, $\{\mu_i\}_{i=1}^r$.

\emph{Coefficients} are defined by ${\bm{\alpha}} = \begin{bmatrix}\alpha_1&\cdots&\alpha_r\end{bmatrix}^T=W^{-1}U_r^*\bm{x}_0$.

We can now reconstruct the approximate dynamics as,
\begin{equation}\label{APPeq:dynamicApproDiscrete}
    \bm{\tilde{x}}_k\approx \Phi D^k {\bm{\alpha}}=\sum_{i=1}^r\alpha_i\mu_i^k\bm{\phi}_i.
\end{equation}

\paragraph{Reconstruction error}
Many applications are satisfied with the above step for recovering the spatial structures in the dynamics. However, for recovering the dynamic with \ac{DMD} another measurement  must be considered. To assess the accuracy, not only the ``moving'' from one sample to the next one should be taken under considerations  but also the dynamic in general. Namely, the criterion should be the summation over the distance between $x_k$ and $\hat{x}_k$. For example, the summation over squared Euclidean distances is resulted in
\begin{equation}\label{eq:errorRec}
    E_{rec}=\sum_{k=0}^M\norm{x_k-\hat{x}_k}^2=\norm{X-\hat{X}}_\mathcal{F}^2
\end{equation}
which is Frobenius norm of the error.

\section{Sparse Representation} \label{appsec:naive}
The main focus should be put on the time profile of the dynamic since the Koopman theory is based on that. In addition, we assume the dynamics induces a family of monotonic time profiles, $\mathcal{D}$, where they differ by their parameters. For example, in linear systems, these functions are exponential, in zero-homogeneous dynamical systems the functions are linear with different slops. 

We assume the typical decay profile is known and we find the nonzero mode for example with the Lasso algorithm \cite{mairal2014sparse}. Then, we remove the not relevant modes and the corresponding atoms in the dictionary. We elaborate the algorithm in Algo. \ref{algo:SR}

\begin{algorithm}[phtb!] \caption{Sparse Representation}
\begin{algorithmic}[1]
		\Inputs{Data sequence $\{{\bm{x}_k}\}_0^{M}$ and decay dictionary $\mathcal{D}$}
		\Initialize{$\mathcal{SR}=\emptyset$}
		\State Find the sparse representation $V$ according to \cite{mairal2014spams}.
		\State Let $\mathcal{I}$ be the set of indices of the atoms in $\mathcal{D}$ sorted (from low to high) according to the norm of the modes (column vectors of $V$).
		\State Remove from $\mathcal{I}$ the indices for which the modes are zeros.
        
		\While{True}
		\State Define $\hat{\mathcal{D}}$ as a new dictionary containing the atoms with indices $\mathcal{I}$.
        \State Compute $\hat{V}=X\hat{\mathcal{D}}^T(\hat{\mathcal{D}}\hat{\mathcal{D}}^T)^{-1}$
		\State Compute the error $\norm{X-\hat{V}\hat{\mathcal{D}}}^2_\mathcal{F}$ 
		\State Add the set $\mathcal{I}$ and its corresponding error to $\mathcal{SR}$
		\State Remove the first index in $\mathcal{I}$.
		\If{$\mathcal{I}$ is empty}
		\State Break
		\EndIf
		\EndWhile
		\State Find in $\mathcal{SR}$ the set of indices $\mathcal{I}$ that yields the minimum error
		\State Define $\hat{\mathcal{D}}$ as a new dictionary containing the atoms with indices $\mathcal{I}$.
        \State Compute $\hat{V}=X\hat{\mathcal{D}}^T(\hat{\mathcal{D}}\hat{\mathcal{D}}^T)^{-1}$
		\Outputs{$$\hat{V},\mathcal{\hat{D}}$$}
    \end{algorithmic}
    \label{algo:SR}
\end{algorithm}

\appendix  
\bibliography{smartPeople}

\begin{thebibliography}{57}
\expandafter\ifx\csname natexlab\endcsname\relax\def\natexlab#1{#1}\fi
\expandafter\ifx\csname selectlanguage\endcsname\relax
  \def\selectlanguage#1{\relax}\fi

\bibitem[\protect\citename{Andreu {et~al.}, }2001]{andreu2001minimizing}
Andreu, Fuensanta, Ballester, Coloma, Caselles, Vicent, and Maz{\'o}n,
  Jos{\'e}~M. 2001.
\newblock Minimizing total variation flow.
\newblock {\em Differential and integral equations}, {\bf 14}(3), 321--360.

\bibitem[\protect\citename{Askham and Kutz, }2018]{askham2018variable}
Askham, Travis, and Kutz, J~Nathan. 2018.
\newblock Variable projection methods for an optimized dynamic mode
  decomposition.
\newblock {\em SIAM Journal on Applied Dynamical Systems}, {\bf 17}(1),
  380--416.

\bibitem[\protect\citename{Azencot {et~al.}, }2019]{azencot2019consistent}
Azencot, Omri, Yin, Wotao, and Bertozzi, Andrea. 2019.
\newblock Consistent dynamic mode decomposition.
\newblock {\em SIAM Journal on Applied Dynamical Systems}, {\bf 18}(3),
  1565--1585.

\bibitem[\protect\citename{Bellettini {et~al.}, }2002]{bellettini2002total}
Bellettini, Giovanni, Caselles, Vicent, and Novaga, Matteo. 2002.
\newblock The total variation flow in RN.
\newblock {\em Journal of Differential Equations}, {\bf 184}(2), 475--525.

\bibitem[\protect\citename{Bollt, }2021]{bollt2021geometric}
Bollt, Erik~M. 2021.
\newblock Geometric Considerations of a Good Dictionary for Koopman Analysis of
  Dynamical Systems: Cardinality,“Primary Eigenfunction,” and Efficient
  Representation.
\newblock {\em Communications in Nonlinear Science and Numerical Simulation},
  105833.

\bibitem[\protect\citename{Brezis, }1973]{brezis1973ope}
Brezis, Haim. 1973.
\newblock {\em Op\'erateurs maximaux monotones et semi-groupes de contractions
  dans les espaces de {H}ilbert}.
\newblock Elsevier.

\bibitem[\protect\citename{Brunton and Kutz, }2019]{brunton2019data}
Brunton, Steven~L, and Kutz, J~Nathan. 2019.
\newblock {\em Data-driven science and engineering: Machine learning, dynamical
  systems, and control}.
\newblock Cambridge University Press.
\newblock Pages  276--320.

\bibitem[\protect\citename{Brunton {et~al.}, }2016]{brunton2016discovering}
Brunton, Steven~L, Proctor, Joshua~L, and Kutz, J~Nathan. 2016.
\newblock Discovering governing equations from data by sparse identification of
  nonlinear dynamical systems.
\newblock {\em Proceedings of the national academy of sciences}, {\bf 113}(15),
  3932--3937.

\bibitem[\protect\citename{Brunton {et~al.}, }2021]{brunton2021modern}
Brunton, Steven~L, Budi{\v{s}}i{\'c}, Marko, Kaiser, Eurika, and Kutz,
  J~Nathan. 2021.
\newblock Modern Koopman theory for dynamical systems.
\newblock {\em arXiv preprint arXiv:2102.12086}.

\bibitem[\protect\citename{Bungert and Burger, }2019]{bungert2019asymptotic}
Bungert, Leon, and Burger, Martin. 2019.
\newblock Asymptotic profiles of nonlinear homogeneous evolution equations of
  gradient flow type.
\newblock {\em Journal of Evolution Equations},  1--32.

\bibitem[\protect\citename{Bungert {et~al.}, }2019a]{bungert2019computing}
Bungert, Leon, Burger, Martin, and Tenbrinck, Daniel. 2019a.
\newblock Computing nonlinear eigenfunctions via gradient flow extinction.
\newblock {Pages  291--302 of:} {\em International Conference on Scale Space
  and Variational Methods in Computer Vision}.
\newblock Springer.

\bibitem[\protect\citename{Bungert {et~al.}, }2019b]{bungert2019nonlinear}
Bungert, Leon, Burger, Martin, Chambolle, Antonin, and Novaga, Matteo. 2019b.
\newblock Nonlinear spectral decompositions by gradient flows of
  one-homogeneous functionals.
\newblock {\em arXiv preprint arXiv:1901.06979}.

\bibitem[\protect\citename{Burger {et~al.}, }2016]{burger2016spectral}
Burger, Martin, Gilboa, Guy, Moeller, Michael, Eckardt, Lina, and Cremers,
  Daniel. 2016.
\newblock Spectral decompositions using one-homogeneous functionals.
\newblock {\em SIAM Journal on Imaging Sciences}, {\bf 9}(3), 1374--1408.

\bibitem[\protect\citename{Chuaqui, }2018]{chuaqui2018general}
Chuaqui, Martin. 2018.
\newblock General criteria for curves to be simple.
\newblock {\em Journal of Mathematical Analysis and Applications}, {\bf
  464}(1), 955--963.

\bibitem[\protect\citename{Cohen and Gilboa, }2018]{cohen:hal-01870019}
Cohen, Ido, and Gilboa, Guy. 2018 (Oct.).
\newblock {\em {Shape Preserving Flows and the p--Laplacian Spectra}}.
\newblock working paper or preprint.

\bibitem[\protect\citename{Cohen and Gilboa, }2020]{cohen2020introducing}
Cohen, Ido, and Gilboa, Guy. 2020.
\newblock Introducing the p-Laplacian spectra.
\newblock {\em Signal Processing}, {\bf 167}, 107281.

\bibitem[\protect\citename{Cohen {et~al.}, }2021a]{cohen2021modes}
Cohen, Ido, Azencot, Omri, Lifshits, Pavel, and Gilboa, Guy. 2021a.
\newblock Modes of homogeneous gradient flows.
\newblock {\em SIAM Journal on Imaging Sciences}, {\bf 14}(3), 913--945.

\bibitem[\protect\citename{Cohen {et~al.}, }2021b]{cohen2021Total}
Cohen, Ido, Berkov, Tom, and Gilboa, Guy. 2021b.
\newblock Total-Variation Mode Decomposition.
\newblock {Pages  52--64 of:} Elmoataz, Abderrahim, Fadili, Jalal, Qu{\'e}au,
  Yvain, Rabin, Julien, and Simon, Lo{\"i}c (eds), {\em Scale Space and
  Variational Methods in Computer Vision}.
\newblock Cham: Springer International Publishing.

\bibitem[\protect\citename{Courant and John, }2012]{courant2012introduction}
Courant, Richard, and John, Fritz. 2012.
\newblock {\em Introduction to calculus and analysis I}.
\newblock Springer Science \& Business Media.

\bibitem[\protect\citename{Dawson {et~al.}, }2016]{dawson2016characterizing}
Dawson, Scott~TM, Hemati, Maziar~S, Williams, Matthew~O, and Rowley,
  Clarence~W. 2016.
\newblock Characterizing and correcting for the effect of sensor noise in the
  dynamic mode decomposition.
\newblock {\em Experiments in Fluids}, {\bf 57}(3), 42.

\bibitem[\protect\citename{Elad, }2010]{elad2010sparse}
Elad, Michael. 2010.
\newblock {\em Sparse and redundant representations: from theory to
  applications in signal and image processing}.
\newblock Springer Science \& Business Media.

\bibitem[\protect\citename{Evangelisti, }2011]{evangelisti2011controllability}
Evangelisti, E. 2011.
\newblock {\em Controllability and Observability: Lectures given at a Summer
  School of the Centro Internazionale Matematico Estivo (CIME) held in
  Pontecchio (Bologna), Italy, July 1-9, 1968}.
\newblock  Vol. 46.
\newblock Springer Science \& Business Media.

\bibitem[\protect\citename{Gavish and Donoho, }2014]{gavish2014optimal}
Gavish, Matan, and Donoho, David~L. 2014.
\newblock The optimal hard threshold for singular values is $4/\sqrt{3}$.
\newblock {\em IEEE Transactions on Information Theory}, {\bf 60}(8),
  5040--5053.

\bibitem[\protect\citename{Gilboa, }2013]{gilboa2013spectral}
Gilboa, Guy. 2013.
\newblock A spectral approach to total variation.
\newblock {Pages  36--47 of:} {\em International Conference on Scale Space and
  Variational Methods in Computer Vision}.
\newblock Springer.

\bibitem[\protect\citename{Gilboa, }2014]{gilboa2014total}
Gilboa, Guy. 2014.
\newblock A total variation spectral framework for scale and texture analysis.
\newblock {\em SIAM journal on Imaging Sciences}, {\bf 7}(4), 1937--1961.

\bibitem[\protect\citename{Gilboa, }2018]{gilboa2018nonlinear}
Gilboa, Guy. 2018.
\newblock {\em Nonlinear Eigenproblems in Image Processing and Computer
  Vision}.
\newblock Springer.

\bibitem[\protect\citename{Gilboa and Osher, }2009]{gilboa2009nonlocal}
Gilboa, Guy, and Osher, Stanley. 2009.
\newblock Nonlocal operators with applications to image processing.
\newblock {\em Multiscale Modeling \& Simulation}, {\bf 7}(3), 1005--1028.

\bibitem[\protect\citename{Hemati {et~al.}, }2017]{hemati2017biasing}
Hemati, Maziar~S, Rowley, Clarence~W, Deem, Eric~A, and Cattafesta, Louis~N.
  2017.
\newblock De-biasing the dynamic mode decomposition for applied Koopman
  spectral analysis of noisy datasets.
\newblock {\em Theoretical and Computational Fluid Dynamics}, {\bf 31}(4),
  349--368.

\bibitem[\protect\citename{Kaiser {et~al.}, }2018]{kaiser2018discovering}
Kaiser, Eurika, Kutz, J~Nathan, and Brunton, Steven~L. 2018.
\newblock Discovering conservation laws from data for control.
\newblock {Pages  6415--6421 of:} {\em 2018 IEEE Conference on Decision and
  Control (CDC)}.
\newblock IEEE.

\bibitem[\protect\citename{Kaiser {et~al.}, }2021]{kaiser2021data}
Kaiser, Eurika, Kutz, J~Nathan, and Brunton, Steven. 2021.
\newblock Data-driven discovery of Koopman eigenfunctions for control.
\newblock {\em Machine Learning: Science and Technology}.

\bibitem[\protect\citename{Katzir, }2017]{Katzir2017Thesis}
Katzir, Oren. 2017 (March).
\newblock {\em On the scale-space of filters and their applications}.
\newblock M.Phil. thesis, Technion — Israel Institute of Technology, Haifa
  3200003.

\bibitem[\protect\citename{Kawahara, }2016]{kawahara2016dynamic}
Kawahara, Yoshinobu. 2016.
\newblock Dynamic mode decomposition with reproducing kernels for Koopman
  spectral analysis.
\newblock {\em Advances in neural information processing systems}, {\bf 29},
  911--919.

\bibitem[\protect\citename{Koopman, }1931]{koopman1931hamiltonian}
Koopman, Bernard~O. 1931.
\newblock Hamiltonian systems and transformation in Hilbert space.
\newblock {\em Proceedings of the national academy of sciences of the united
  states of america}, {\bf 17}(5), 315.

\bibitem[\protect\citename{Korda and Mezi{\'c}, }2018]{korda2018linear}
Korda, Milan, and Mezi{\'c}, Igor. 2018.
\newblock Linear predictors for nonlinear dynamical systems: Koopman operator
  meets model predictive control.
\newblock {\em Automatica}, {\bf 93}, 149--160.

\bibitem[\protect\citename{Kutz {et~al.}, }2016a]{kutz2016dynamic}
Kutz, J~Nathan, Brunton, Steven~L, Brunton, Bingni~W, and Proctor, Joshua~L.
  2016a.
\newblock {\em Dynamic mode decomposition: data-driven modeling of complex
  systems}.
\newblock SIAM.

\bibitem[\protect\citename{Kutz {et~al.}, }2016b]{kutz2016koopman}
Kutz, J~Nathan, Proctor, Joshua~L, and Brunton, Steven~L. 2016b.
\newblock Koopman theory for partial differential equations.
\newblock {\em arXiv preprint arXiv:1607.07076}.

\bibitem[\protect\citename{Langley {et~al.}, }1981]{langley1981bacon}
Langley, Pat, Bradshaw, Gary~L, and Simon, Herbert~A. 1981.
\newblock BACON. 5: The discovery of conservation laws.
\newblock {Pages  121--126 of:} {\em IJCAI},  vol. 81.
\newblock Citeseer.

\bibitem[\protect\citename{Li {et~al.}, }2017]{li2017extended}
Li, Qianxiao, Dietrich, Felix, Bollt, Erik~M, and Kevrekidis, Ioannis~G. 2017.
\newblock Extended dynamic mode decomposition with dictionary learning: A
  data-driven adaptive spectral decomposition of the Koopman operator.
\newblock {\em Chaos: An Interdisciplinary Journal of Nonlinear Science}, {\bf
  27}(10), 103111.

\bibitem[\protect\citename{Lu and Tartakovsky, }2020]{lu2020prediction}
Lu, Hannah, and Tartakovsky, Daniel~M. 2020.
\newblock Prediction accuracy of dynamic mode decomposition.
\newblock {\em SIAM Journal on Scientific Computing}, {\bf 42}(3),
  A1639--A1662.

\bibitem[\protect\citename{Mairal {et~al.}, }2014a]{mairal2014spams}
Mairal, Julien, Bach, F, Ponce, J, Sapiro, G, Jenatton, R, and Obozinski, G.
  2014a.
\newblock Spams: A sparse modeling software, v2. 6.
\newblock {\em URL http://spams-devel. gforge. inria. fr/downloads. html}.

\bibitem[\protect\citename{Mairal {et~al.}, }2014b]{mairal2014sparse}
Mairal, Julien, Bach, Francis, Ponce, Jean, {et~al.} 2014b.
\newblock Sparse Modeling for Image and Vision Processing.
\newblock {\em Foundations and Trends{\textregistered} in Computer Graphics and
  Vision}, {\bf 8}(2-3), 85--283.

\bibitem[\protect\citename{Mauroy, }2021]{mauroy2021koopman}
Mauroy, Alexandre. 2021.
\newblock Koopman Operator Theory for Infinite-Dimensional Systems: Extended
  Dynamic Mode Decomposition and Identification of Nonlinear PDEs.
\newblock {\em arXiv preprint arXiv:2103.12458}.

\bibitem[\protect\citename{Mauroy {et~al.}, }2020]{mauroy2020koopman}
Mauroy, Alexandre, Susuki, Y, and Mezi{\'c}, I. 2020.
\newblock {\em The Koopman Operator in Systems and Control}.
\newblock Springer.

\bibitem[\protect\citename{Mezi{\'c}, }2005]{mezic2005spectral}
Mezi{\'c}, Igor. 2005.
\newblock Spectral properties of dynamical systems, model reduction and
  decompositions.
\newblock {\em Nonlinear Dynamics}, {\bf 41}(1-3), 309--325.

\bibitem[\protect\citename{Nakao and Mezi{\'c}, }2020]{nakao2020spectral}
Nakao, Hiroya, and Mezi{\'c}, Igor. 2020.
\newblock Spectral analysis of the Koopman operator for partial differential
  equations.
\newblock {\em Chaos: An Interdisciplinary Journal of Nonlinear Science}, {\bf
  30}(11), 113131.

\bibitem[\protect\citename{Nathan~Kutz {et~al.}, }2018]{nathan2018applied}
Nathan~Kutz, J, Proctor, Joshua~L, and Brunton, Steven~L. 2018.
\newblock Applied Koopman theory for partial differential equations and
  data-driven modeling of spatio-temporal systems.
\newblock {\em Complexity}, {\bf 2018}.

\bibitem[\protect\citename{Otto and Rowley, }2021]{otto2021koopman}
Otto, Samuel~E, and Rowley, Clarence~W. 2021.
\newblock Koopman operators for estimation and control of dynamical systems.
\newblock {\em Annual Review of Control, Robotics, and Autonomous Systems},
  {\bf 4}.

\bibitem[\protect\citename{Pan {et~al.}, }2021]{pan2021sparsity}
Pan, Shaowu, Arnold-Medabalimi, Nicholas, and Duraisamy, Karthik. 2021.
\newblock Sparsity-promoting algorithms for the discovery of informative
  Koopman-invariant subspaces.
\newblock {\em Journal of Fluid Mechanics}, {\bf 917}.

\bibitem[\protect\citename{Rudy {et~al.}, }2017]{rudy2017data}
Rudy, Samuel~H, Brunton, Steven~L, Proctor, Joshua~L, and Kutz, J~Nathan. 2017.
\newblock Data-driven discovery of partial differential equations.
\newblock {\em Science Advances}, {\bf 3}(4), e1602614.

\bibitem[\protect\citename{Schmid, }2010]{schmid2010dynamic}
Schmid, Peter~J. 2010.
\newblock Dynamic mode decomposition of numerical and experimental data.
\newblock {\em Journal of fluid mechanics}, {\bf 656}, 5--28.

\bibitem[\protect\citename{Schmidt and Lipson, }2009]{schmidt2009distilling}
Schmidt, Michael, and Lipson, Hod. 2009.
\newblock Distilling free-form natural laws from experimental data.
\newblock {\em science}, {\bf 324}(5923), 81--85.

\bibitem[\protect\citename{Tu {et~al.}, }2013]{tu2013dynamic}
Tu, Jonathan~H, Rowley, Clarence~W, Luchtenburg, Dirk~M, Brunton, Steven~L, and
  Kutz, J~Nathan. 2013.
\newblock On dynamic mode decomposition: Theory and applications.
\newblock {\em arXiv preprint arXiv:1312.0041}.

\bibitem[\protect\citename{Valmorbida and Anderson,
  }2017]{valmorbida2017region}
Valmorbida, Giorgio, and Anderson, James. 2017.
\newblock Region of attraction estimation using invariant sets and rational
  Lyapunov functions.
\newblock {\em Automatica}, {\bf 75}, 37--45.

\bibitem[\protect\citename{Venturi and Dektor, }2021]{venturi2021spectral}
Venturi, Daniele, and Dektor, Alec. 2021.
\newblock Spectral methods for nonlinear functionals and functional
  differential equations.
\newblock {\em Research in the Mathematical Sciences}, {\bf 8}(2), 1--39.

\bibitem[\protect\citename{Williams {et~al.}, }2015a]{williams2015data2}
Williams, Matthew~O, Kevrekidis, Ioannis~G, and Rowley, Clarence~W. 2015a.
\newblock A data--driven approximation of the koopman operator: Extending
  dynamic mode decomposition.
\newblock {\em Journal of Nonlinear Science}, {\bf 25}(6), 1307--1346.

\bibitem[\protect\citename{Williams {et~al.}, }2015b]{williams2015data}
Williams, Matthew~O, Rowley, Clarence~W, Mezi{\'c}, Igor, and Kevrekidis,
  Ioannis~G. 2015b.
\newblock Data fusion via intrinsic dynamic variables: An application of
  data-driven Koopman spectral analysis.
\newblock {\em EPL (Europhysics Letters)}, {\bf 109}(4), 40007.

\bibitem[\protect\citename{Williams {et~al.}, }2016]{williams2016extending}
Williams, Matthew~O, Hemati, Maziar~S, Dawson, Scott~TM, Kevrekidis, Ioannis~G,
  and Rowley, Clarence~W. 2016.
\newblock Extending data-driven Koopman analysis to actuated systems.
\newblock {\em IFAC-PapersOnLine}, {\bf 49}(18), 704--709.

\end{thebibliography}
\section*{Acknowledgements}
We would like to thank Prof. Gershon Wolansky from
Department of Mathematics, Technion and Dr. Eli Appelboim from Electrical and Computer Engineering Department, and Dan Glaubach for stimulating discussions.



\end{document}